\documentclass[pdflatex,sn-mathphys-num]{sn-jnl}

\usepackage{graphicx}%
\usepackage{multirow}%
\usepackage{amsmath,amssymb,amsfonts}%
\usepackage{amsthm}%
\usepackage{mathrsfs}%
\usepackage[title]{appendix}%
\usepackage{xcolor}%
\usepackage{textcomp}%
\usepackage{manyfoot}%
\usepackage{booktabs}%
\usepackage{algorithm}%
\usepackage{algorithmicx}%
\usepackage{algpseudocode}%
\usepackage{listings}%
\usepackage{subcaption}
\usepackage{array}
\usepackage{adjustbox}
\usepackage{multirow}
\usepackage{mathdots}
\usepackage{stfloats}
\usepackage{url}
\usepackage{verbatim}
\usepackage{algpseudocode}
\usepackage{rotating} 
\setcitestyle{authoryear,round,aysep={ }}

\theoremstyle{thmstyleone}%
\theoremstyle{thmstyletwo}%

\theoremstyle{thmstylethree}%

\begin{document}

\title[Article Title]{A bi-level priority sorting framework for flexible AGV service scheduling in smart warehouses}


\author*[1]{\fnm{Xiaozhu} \sur{Sun}}\email{xiaozhu.sun@torontomu.ca}

\author[2]{\fnm{Bilal} \sur{Farooq}}\email{bilal.farooq@torontomu.ca}


\affil*[1]{\orgdiv{Laboratory of Innovations in Transportation (LiTrans)}, \orgname{Toronto Metropolitan University}, \orgaddress{\street{350 Victoria St}, \city{Toronto}, \postcode{M5B 2K3}, \state{Ontario}, \country{Canada}}}

\affil[2]{\orgdiv{Laboratory of Innovations in Transportation (LiTrans)}, \orgname{Toronto Metropolitan University}, \orgaddress{\street{350 Victoria St}, \city{Toronto}, \postcode{M5B 2K3}, \state{Ontario}, \country{Canada}}}



\abstract{This paper proposes a bi-level optimization framework to coordinate Automated Guided Vehicle (AGV) flexible operations in smart independent warehouses, addressing the critical challenge of balancing high-throughput order fulfillment with stringent cost control. The framework is designed to simultaneously optimize flexible customer service level, system cost, and operational efficiency. The first level dynamically adjusts real-time scheduling parameters, such as order commitment times and delay tolerance, based on predefined customer priority categories. The second level performs real-time routing optimization for each AGV by identifying the shortest feasible paths while avoiding conflicts. For complex multi-capacity package picking tasks, two heuristic rules, priority, deadline, with shortest path (PDSP) and delay cost with shortest path (DCSP), are applied to multi-capacity package picking tasks and further training is carried out using the reinforcement learning algorithm of A* guided deep Q-learning (AGDQN). Comprehensive simulation experiments, conducted across diverse warehouse layouts and order demand patterns, demonstrate that the proposed framework equipped with both heuristic rules consistently reduces average order delay and total system costs by over 50\% during peak demand periods. This is achieved while maintaining a service level above 90\% and maximizing AGV utilization. The method also exhibits superior flexibility and sustained efficiency under normal and fluctuating demand scenarios. Additional ablation studies confirm that the proposed priority sorting mechanism delivers robust performance advantages when tested with various other reinforcement learning baselines.}




\maketitle

\section{Introduction}\label{sec1}

Logistics automation is an important link in the supply chain, ensuring the efficient circulation of commodities \citep{Yan2022}. The rapid expansion of e-commerce, express delivery, smart factories, and smart cities has propelled intelligent logistics to the forefront of industrial and academic research fields. Within this domain, the integration of smart warehousing systems with automated guided vehicles (AGVs) represents a transformative advancement, significantly enhancing operational efficiency and effectiveness \citep{Fazlollahtabar2015}. As autonomous material-handling devices, AGVs are extensively deployed in manufacturing\citep{LinJamesT.2019Soaf}, automated storage and retrieval systems (AS/RS)\citep{CanoJoseAlejandro2023Stpr}, and port terminals\citep{XieTianhao2025AsaU}. In automated sorting operations, AGVs offer superior endurance and consistent performance compared with manual or semi-automated alternatives. The effective scheduling of AGV fleets is paramount, enabling continuous operation and ensuring timely order fulfillment, thereby elevating overall customer service levels \citep{liu2019multi}.

The optimization of warehouse operations is economically crucial, with order picking alone accounting for approximately 55\% of the total operational costs \citep{Larco02112017}. Research on AGV-based systems has evolved in both macro- and micro-dimensions. At the macro (system) level, studies have focused on integrated challenges such as multi-AGV task scheduling, cooperative path planning, and optimal fleet sizing to maximize throughput by \cite{SunYangjun2025Bfmp} and \cite{ChaikovskaiaMari2025Ooaf}. At the micro-level, investigations delve into real-time collision avoidance, energy consumption minimization, and motion control. In flexible manufacturing and logistics contexts, key performance objectives consistently include minimizing makespan, optimizing delivery routes, reducing work-in-process inventory, and ensuring timely material delivery by \cite{Zhang2023} and \cite{riazi2020energy}.

This study is motivated by a pressing real-world challenge in e-commerce logistics: the need to dynamically balance system efficiency with differentiated customer service. Modern sales and inventory strategies often necessitate the assignment of dual priority levels to orders to cater to diverse customer segments (e.g., premium vs. standard) by \cite{Vicil2021} and \cite{uzunoglu2020continuous}. This creates a complex operational dilemma in which a homogeneous ``first-come, first-served'' AGV scheduling strategy becomes suboptimal. The core problem transcends simple path planning and involves a bi-level decision-making process: first, dynamically adjusting scheduling priorities based on real-time order attributes and customer categories; second, optimizing the physical routing of a fleet of multi-capacity AGVs under those priorities to minimize systemic costs and delays.

To address this gap, we propose a novel bilevel optimization framework that enhances the classical dial-a-ride problem (DARP) model by \cite{Hyytia2012} for smart warehouse sorting. Our framework explicitly models the trade-off between global system efficiency (e.g., total cost and AGV utilization) and customer-oriented metrics (e.g., service level and delay). A weighted cost model was employed to integrate these competing objectives. At its core is a bi-level priority sorting framework: the first level handles the external optimization of order scheduling parameters (such as deadline and delay tolerance) per customer category; the second level conducts the internal optimization of AGV routing by calculating the shortest feasible paths that satisfy the high-level priorities. The main contributions are following:\\

\begin{enumerate}
    \item \emph{Integration of dynamic order priorities and delay constraints into multi-load AGV scheduling rules:} We propose two heuristic strategies, DCSP (delay cost and shortest path) and PDSP (priority, deadline, and shortest path), which explicitly incorporate real-time order priority and delivery delay considerations into the dispatching logic of multi-load AGVs. These rules jointly optimize energy consumption and customer satisfaction, providing a practical and scalable decision-making flexible framework for dynamic warehouse environments.
    \item \emph{Dual-metric evaluation framework based on operational workload and service Level:} This study introduces a systematic evaluation approach that links operational workload with service-level performance. For the first time, this framework was applied to continuous multi‑node picking using automated guided vehicles (AGVs), enabling comprehensive sensitivity analysis, scalability testing, and system state assessment. The results offer actionable insights for the scheduling model design in cloud warehouses, complex task environments, and distribution centers.
    \item \emph{Extensive experimental validation of the bilevel optimization architecture:} Under consistent experimental settings, the proposed method was compared with several state-of-the-art reinforcement learning algorithms. The results demonstrate its superior adaptability, sensitivity, and robustness, confirming that the coordinated two-level optimization effectively enhances the overall system performance in terms of time, cost, and energy efficiency.
\end{enumerate}

The remainder of this paper is organized as follows. Section \ref{sec:LR} reviews the relevant literature. Section \ref{problemdescription} describes the application scenarios, underlying assumptions and problem formulation. Section \ref{model} presents the proposed mathematical model and optimization methods. Section \ref{case} presents the case studies and comparative performance analyses. Finally, Section \ref{conclusion} concludes the paper and discusses the future research directions.

\section{Related Works}
\label{sec:LR}
In the existing literature, AGVs have been extensively studied in manufacturing workshops, focusing on order sorting, priority setting, delivery timeliness, and dynamic programming. However, most research on AGV applications in warehouses emphasizes minimizing path costs, addressing single-task operations, managing shelf transfer vehicles, and solving operational challenges such as deadlocks and congestion. These studies often assumed static scheduling frameworks. This study draws on three major categories of relevant research: recent developments in tactical optimization methods, system analysis and evaluation techniques, including commonly used performance indicators, and classic and contemporary sorting strategies in smart logistics environments. The limitations of existing multi-priority sorting methods are also highlighted, with a comparative analysis of prior models to identify research gaps.
\subsection{Tactical Optimization}
Strategic issues in AGV path optimization often involve vehicle design, facility layout, and route selection, which are typically addressed using multi-objective models. 
\cite{Li2019} proposed a three-objective model using harmony search to reduce AGV travel time. 
\cite{Zhang2021} developed an energy-aware optimization model that considers turning energy consumption and evaluated it using particle swarm optimization (PSO). 
\cite{Zhang2024} proposed a stage-based dynamic prioritization and rerouting approach to reduce congestion and delays in multi-AGV systems. 
\cite{Nishi2011} applied mixed-integer linear programming (MILP) to decouple high-level production scheduling from low-level AGV routing. 
\cite{Wang2020} employed queuing theory within a two-stage robotic mobile fulfillment system (RMFS) framework to analyze system performance.

Once target nodes are determined, enhancing AGV responsiveness and transport efficiency often focuses on shortest-path algorithms. 
\cite{Guo2023} enhanced the traditional A* algorithm by reducing sharp turns in single-AGV paths and minimizing intersection conflicts among multiple AGVs through the introduction of penalty values for overlapping path segments. 
\cite{Zhang2024} applied the A* algorithm to evaluate different scheduling modes involving inter-AGV transfers. 
\cite{XieTianhao2025AsaU} proposed a GAHA algorithm for dynamic and conflict-free AGV scheduling in container terminals to improve efficiency. 
\cite{Jiang2024} addressed AGV routing problems with simultaneous pickup and delivery under mixed time-window constraints. 
\cite{liu2019multi} proposed a multi-objective AGV task scheduling model integrated with adaptive genetic algorithms (AGA and MAGA), with scheduling results visualized using Gantt charts based on a practical sorting layout. 
These methods collectively support greater AGV flexibility in smart logistics systems.

\cite{durst2023multi} developed an adaptive deep reinforcement learning (DRL)-based real-time AGV scheduling approach to optimize multiple efficiency-related objectives in a job shop environment, employing a linear reward function and proximal policy optimization (PPO) with a deep neural network (DNN) for policy iteration. 
Similarly, \cite{LiuChien-Liang2020ADRL} proposed a DRL framework integrating actor and critic networks with convolutional and fully connected layers, enabling agents to learn appropriate behaviors across diverse operational scenarios while estimating state values through a parallel training strategy that combines asynchronous updates with deep deterministic policy gradient (DDPG). 
In their approach, the reward function was manually designed and required adjustment across different environments during the optimization process.

\subsection{System Assessment}\label{LR2} 

AGV performance in smart warehousing systems must be evaluated differently from traditional vehicle metrics, encompassing AGV-specific performance degradation, inventory impacts, and responsiveness to orders and service sequences. Inventory optimization has yielded significant results in this field. For instance, \cite{Vicil2021} introduced two ordering policy models, one of which assumed zero tolerance for waiting among high-priority customers and used continuous-time Markov chain (CTMC) modeling for approximation. \cite{Kocer2020} investigated a lost-sales (s,Q) inventory system serving two customer classes. By comparing equal and graded policies, they concluded that prioritizing high-priority customers yields lower supplier costs and better service levels than the latter. \cite{chen2010flexible} emphasized that selecting appropriate priority strategies is complex, requiring trade-offs among handled volume, service level, and operational costs. 

\cite{Gupta202} modeled dual-warehouse systems with delayed payments; \cite{Goh2014} explored inventory classifications and system costs. Energy and routing efficiency are also important, while \cite{Qiu2015} modeled energy use by weight and distance in high-bay settings. 

Delay-related costs are increasingly relevant. \cite{Gao2020} incorporated soft time window constraints and quantified waiting and delay costs across four urban scenarios. \cite{Arnold2013} examined Maersk Line’s use of an economic framework to accelerate value delivery, analyzing delay costs based on benefits, urgency, and value of information. \cite{Anderson2012} formally introduced the concept of cost of delay (CoD) within Kanban-based product development, arguing that understanding CoD is essential for optimizing work prioritization and customer value delivery. Despite the relevance of these concepts, few studies have combined them with AGV systems.

\subsection{Priority Evaluation}\label{LR3} 
To enhance the flexibility and responsiveness of logistics systems, scholars have differentiated between internal and external priorities. 
\cite{olsson2008external} distinguished internal priorities such as improving process efficiency, optimizing resource allocation, and minimizing costs from external priorities, which include customer requirements, market trends, and delivery deadlines. 
\cite{arshinina2019strategic} set internal priorities serve the long-term operational goals of the organization, while external priorities focus on satisfying immediate customer needs and enhancing market competitiveness.

In AGV systems, the first-come-first-served (FCFS) rule ~\citep{alma992404284905151} is a widely applied classical strategy. 
This rule is used to determine the order in which vehicles are dispatched based on the sequence in which requests arrive, ensuring fairness, but often at the expense of overall efficiency. 
In practice, this rule is the most prevalent queuing strategy in AGV dispatching systems. 
However, many studies have explored alternative scheduling rules to further refine vehicle coordination. 
\cite{Goldberg1977} applied the earliest due date (EDD) rule to manage queues, aiming to reduce delay penalties by prioritizing tasks with the earliest possible deadlines. 
\cite{jackson1955} designed scheduling mechanisms that minimize job tardiness, forming the basis of many modern dispatching approaches. 
~\cite{Smith1956} proposed the shortest processing time (SPT) rule, which prioritizes tasks requiring the least processing time, effectively reducing average waiting time and increasing system throughput.

Beyond these heuristic rules, priority modeling based on cost functions has emerged in queuing and service systems research. 
\cite{Maglaras2018} examined packet delays in single server queues by modeling priority through delay cost functions, enabling the system to prioritize tasks based on the economic impact of delay. 
Similarly, \cite{Amiri1998} introduced external priority attributes into transportation planning, including time costs, expected workload, and backup service requirements. 
Their model allows for more nuanced scheduling decisions that better align with customer and operational requirements.

\subsection{Remarks on Existing Literature} 
Despite significant advances, several critical gaps remain in AGV optimization research. Dynamic scheduling methods in warehouse logistics are underdeveloped compared to those in flexible manufacturing systems. Although traditional transportation and smart factory designs offer some integrated optimization methods, innovative approaches for dynamic multi-capacity, multi-task AGV fleet scheduling in logistics environments are limited. Most existing studies focus on individual AGV scheduling or fleet size reduction, emphasizing algorithmic improvements for the shortest path search or vehicle count minimization.

However, they often overlook external strategic and marginal factors such as customer satisfaction, inventory economics, and priority differentiation. Moreover, there is a lack of research addressing simultaneous improvements in order responsiveness, travel time reduction, and energy efficiency for both loaded and empty states. Moreover, the prevalent assumption of a single-order capacity per AGV does not reflect real-world operations, where AGVs frequently carry multiple orders, significantly boosting operational efficiency.

\section{Problem Description}\label{problemdescription} 
\begin{figure}[h]
\centering
\includegraphics[width=0.9\textwidth]{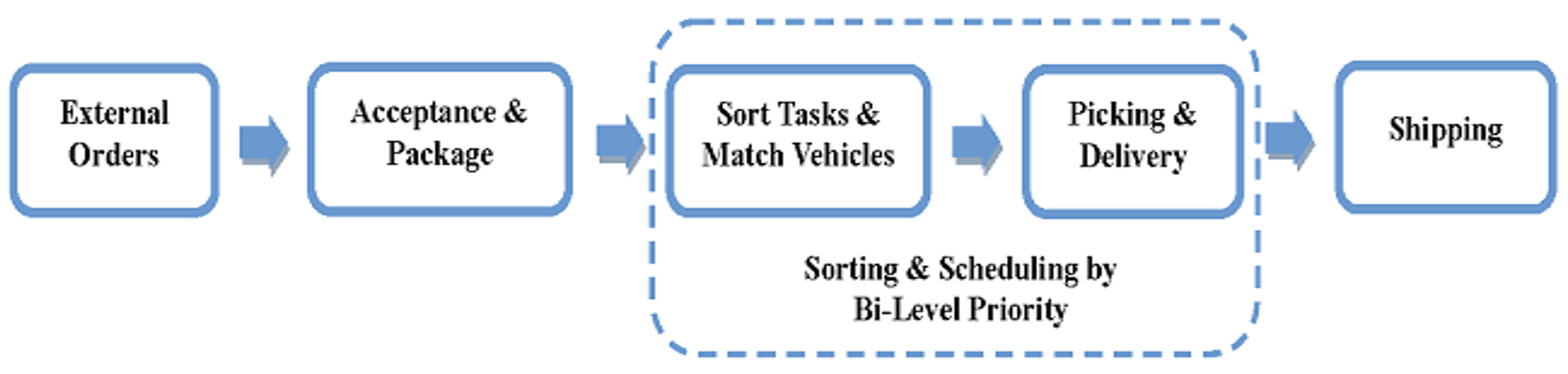}
\caption{\textbf {Delivery processes.}}
\label{fig:processes}
\end{figure}

\begin{figure}[h]
\centering
\includegraphics[width=0.9\textwidth]{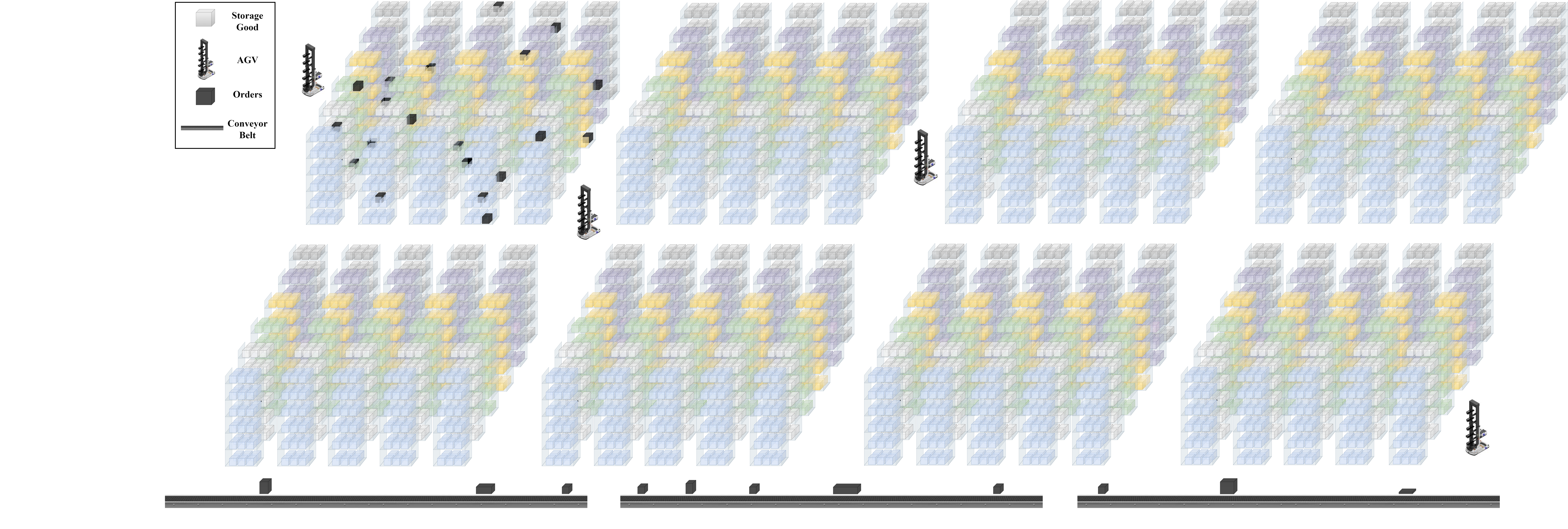}
\caption{\textbf{ Considered partial layout of the warehouse.}}
\label{fig:layout}
\end{figure}

In traditional warehouses, operations follow a sequence of picking, packing, staging, scheduling, and dispatching~\citep{kim2020improving}, whereas in smart warehouses, AS/RS fetch items based on customer orders and deliver them via conveyors to the packaging area~\citep{liu2019multi}. Packaged items are then conveyed to the sorting zone, where AGVs handle the sorting and delivery, as illustrated in Figure \ref{fig:processes}, with the system layout shown in Figure \ref{fig:layout}. Each 24-hour cycle begins with a static planning phase that allocates tasks and schedules any remaining orders from the previous cycle, followed by dynamic planning as new orders and retrieval tasks emerge throughout the day. AGVs continuously update their tasks in real time, retrieving assignments in the delivery zone, anticipating subsequent tasks, and automatically starting new cycles after completing each operational loop.

In a grid-based warehouse map, each delivery location is assigned specific coordinates, and lines denote navigable AGV paths. To minimize backorders and lost sales by~\cite{Vicil2021} and \cite{uzunoglu2020continuous}, the system incorporates customer priority classes, giving higher weights and stricter delay penalties to high-priority customers, whereas lower-priority customers receive lighter penalties. Because the number of AGVs is typically far fewer than the volume of items awaiting transport, queuing and delays are unavoidable, and optimal routing must adapt dynamically to fleet size. The model first considers a single AGV handling all orders, recalculating and reordering tasks as new orders arrive, optimal tasks emerge, or scheduled orders are canceled, and sometimes requires mid-task redirection. It is then extended to multi-AGV coordination, enabling the global optimization of task allocation and routing.

The following assumptions were made to facilitate the path optimization of AGVs in the warehouse environment:
\begin{itemize}

    \item Each AGV carries up to four units, adjustable as needed.

    \item AGVs move at constant speed for safety and predictability.

    \item Loading/unloading times are excluded from transport time.

    \item Corridors allow only one AGV; conflicts are resolved by priority rules.

    \item An AGV must finish its current loop before new tasks.
    
\end{itemize}

\section{Methodology}\label{model} 
This study extends the classical dial-a-ride problem (DARP) by~\cite{Hyytia2012}, traditionally used in demand-responsive transport on urban road networks, to AGV-based warehouse logistics. Both systems involve capacity-constrained and time-sensitive pickup and delivery tasks. Although their evaluation criteria differ slightly, the core objective remains minimizing operational costs while maintaining high service quality. In warehouse micro-transportation, AGVs must navigate dynamic conditions to meet similar efficiency and service requirements. To adapt the DARP to AGV-based logistics, we introduce two key performance metrics: system effort and time-varying level of service (LoS). The system effort accounts for the AGV travel distance, idle time, and energy consumption, reflecting the overall operational efficiency. In contrast, the time-varying LoS evaluates how effectively AGVs meet demand under changing conditions, such as customer wait times, delivery punctuality, and throughput. These metrics enable a comprehensive evaluation of the performance under dynamic logistics scenarios.

\begin{figure}[t!]
\centering
\includegraphics[width=0.8\textwidth]{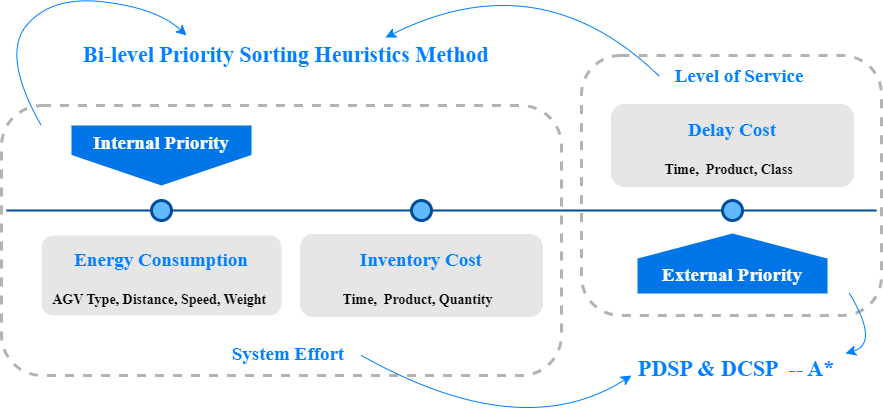}
\caption{\textbf{Concept of the proposed methodology.}}
\label{fig:method}
\end{figure}

A bi-level optimization framework was developed to integrate these metrics into AGV scheduling. The internal level prioritizes task assignments based on distance and inventory fluctuations, whereas the external level reflects service-driven factors such as customer priority, time sensitivity, and delay tolerance. Two model-based sorting mechanisms guide AGV coordination for order-picking and delivery. Depending on the system constraints, minimizing the workload and maximizing the LoS may require trade-offs or produce mutually beneficial outcomes, as illustrated in Figure ~\ref{fig:method}.

This study models AGV routing as a many-to-one ride-matching problem, in which one AGV can carry multiple orders, but each order is assigned to only one AGV. The orders include spatial coordinates, service attributes, priority levels, and destinations. Internal and external priorities guide the allocation to minimize the total system cost, which is evaluated using an enhanced DARP framework. Based on the dynamic DARP theory, the model targets two objectives: minimizing AGV energy consumption (system effort) and reducing the total system time, including waiting and delivery delays (level of service, LoS). As these objectives may conflict, the model minimizes the weighted sum to balance efficiency and service quality.

The detailed notations and related definitions are in Table \ref{tab:notations}, as well as the mathematical formulas, are shown below.

\begin{table}[h]
\footnotesize
\caption{Notations and definitions}
\label{tab:notations}
\begin{tabular*}{\textwidth}{@{\extracolsep\fill}lp{0.08\textwidth}p{0.28\textwidth}p{0.08\textwidth}p{0.28\textwidth}}
\toprule
\textbf{Type} 
& \multicolumn{2}{c}{\textbf{Notation 1}} 
& \multicolumn{2}{c}{\textbf{Notation 2}} \\
\cmidrule(lr){2-3}\cmidrule(lr){4-5}
& \textbf{Symbol} & \textbf{Definition}
& \textbf{Symbol} & \textbf{Definition} \\
\midrule

\multirow{3}{*}{\textsc{Indices}}
& $i$ & Order index, $i=1,\ldots,N$
& $j$ & Load index, $j=0,\ldots,M$ \\
& $cl$ & Priority class, $cl\in\{A,B,C,D\}$
& $k$ & AGV index, $k=1,\ldots,K$ \\
& $t$ & Time window index
& $(i,j)$ & Arc from node $i$ to $j$ \\

\midrule
\multirow{6}{*}{\textsc{Sets}}
& $P_i$ & Orders in one loop
& $p_j$ & Tasks in one trip \\
& $O_i$ & Origin nodes of tasks
& $D_i$ & Destination nodes of tasks \\
& $C_{cl}^i$ & Orders with class $cl$
& $A_k$ & Attributes of AGV $k$ \\
& $\mathcal{O}_t$ & Active orders at time $t$
& $\mathcal{N}$ & Set of warehouse nodes \\
& $\mathcal{K}$ & Set of AGVs
& $\mathcal{A}$ & Set of network arcs \\
& $g_i$ & Task weight (kg)
&  &  \\

\midrule
\multirow{8}{*}{\textsc{Parameters}}
& $N$ & Number of orders
& $M$ & AGV capacity \\
& $K$ & Number of AGVs
& $C$ & Number of priority classes \\
& $G_0$ & AGV self-weight (kg)
& $R$ & Conveyor coordinates \\
& $T_B$ & Battery tolerance time (h)
& $v$ & Velocity (m/s) \\
& $EP$ & Electricity price
& $\delta_1,\delta_2$ & Power coefficients \\
& $C_O$ & Waiting cost
& $C_{cl}^p$ & Delay cost rate \\
& $\beta$ & Penalty factor
& $d_{ij}$ & Distance between nodes $i$ and $j$ (m) \\
& $n_{ij}^{\max}$ & Arc capacity
& $T_{ddi}$ & Order deadline (s) \\

\midrule
\multirow{10}{*}{\textsc{Variables}}
& $t_{oi}$ & Order arrival time
& $t_{ai}$ & Arrival time at pickup \\
& $t_{di}$ & Departure time
& $P_{o i}$ & Power consumption of order $i$ \\
& $P_{unit}^k$ & Power rate of AGV $k$
& $T_{ri}$ & Running time \\
& $GL$ & Full-load weight
& $w_t$ & Dynamic weight \\
& $T_i$ & Total transport time
& $T_{wi}$ & Waiting time \\
& $T_{ui}$ & Uniform waiting baseline
& $d_i$ & Transport distance \\
& $x_{ijk}$ & Routing decision variable
& $y_{ik}$ & Assignment variable \\
& $s_i$ & Service start time
& $u_i$ & Subtour elimination variable \\
& $E_k$ & Energy consumption of AGV $k$
& $C_{Power}$ & Power cost \\
& $C_{Time}$ & Time-related cost
& $G_k$ & Carried weight of AGV $k$ \\

\bottomrule
\end{tabular*}
\end{table}

\subsection{Objective Function}\label{function} 
In each time window, the system receives a batch of new orders with assigned priority levels, prompting the AGVs to re-evaluate and optimize the dispatch routes accordingly. Traditional objective functions often fall short of meeting multi-class customer demands. To overcome this, a dynamic weighting mechanism in Equation \ref{utilityMLM} is introduced, assigning higher weights to urgent orders or applying queue-time penalties. When the waiting times exceed the set thresholds, the penalty coefficients are adjusted based on the customer priority classes.

\begin{equation}
\small
\min  F =  w_t\times{ E[Power]} + (1- w_t)\times{ E[Time]}
\label{utilityMLM}
\end{equation}

where $ E[Power]$ represents the expected monetary cost of the AGV energy consumption, which is influenced by the vehicle velocity, transport distance, and payload weight. $ E[Time]$ denotes the expected cost associated with the system time, including travel, waiting, and delay durations. The parameter $ w_t$ adjusts the trade-off between the energy and time costs. The system dynamically tunes this weight in response to operational factors such as order priority, customer importance, and upstream and downstream supply chain impacts.

\subsubsection{Enhanced Power Function}\label{power} 
The power consumption model in Equation \ref{2}  is designed with flexibility. However, for this study, we assume that AGVs operate at a constant speed while transferring electrical loads, consistent with large-capacity transport systems~\citep{Qiu2015}. As the AGVs were assumed to operate continuously, a motion coefficient was applied for AGVs carrying a payload based on prior work~\citep{Mei2006}.

\begin{equation}
\small
P_{unit}^k =  \delta_1 \times  G_k +  \delta_2
\label{2}
\end{equation}
where the total weight $G_k$ includes AGV tare weight and dynamic payload:

\begin{equation}
G_k = K\times{ G_0 + \sum_{i=1}^{n}  g_i \times y_{ik} + \sum_{(i,j)\in A} q_{ijk}}
\label{total_weight}
\end{equation}

Total power consumption over distance $L_k$ traveled by AGV $k$:

\begin{equation}
Po_k = P_{unit}^k \times L_k
\label{total_power}
\end{equation}

\begin{equation}
L_k = \sum_{(i,j) \in A} d_{ij} \times x_{ijk}
\label{total_distance}
\end{equation}

where $d_{ij}$ is the distance between nodes $i$ and $j$.


\subsubsection{Comprehensive Time Function}\label{time} 
The system defines a uniform waiting time for all orders but with different thresholds for each class. When the waiting time for an order exceeds the threshold, the excess time is classified as a delay. In Equations \ref{6} and \ref{7}, $T(waiting)$ and $T(delay)$ denote the cumulative waiting and delay times for all orders in a given time window.

\begin{equation}
C_{Time} = C_{waiting} + \beta \times C_{delay}
\label{time_function}
\end{equation}

\begin{equation}
C_{Time} = \sum_{i=1}^{n} \left[ f_{C_O}(T_{wi}) + \beta \times f_{C_{cl}^i}(T_{delay}) \right]
\label{time_function_complete}
\end{equation}

\begin{equation}
\small
T(waiting) =  \sum_{i=1}^{n} T_{wi}
\label{6}
\end{equation}
where $T_{wi} = \max(0, s_i - r_i)$ is the waiting time for order $i$, $s_i$ is the service start time and $r_i$ is the request time.
\begin{equation}
\small
T(delay) = \sum_{i=1}^{n}  |T_{ddi} - T_{wi}|
\label{7}
\end{equation}

\subsubsection{Conversion Function}\label{system} 
To enable cost normalization, both the power and time metrics were converted into monetary values. The electricity cost was calculated using the unit industrial power rate, whereas the time-related costs were derived based on order-specific penalties. Similarly, we obtain unit costs for different time types and calculate the lost costs. $C_ O$ is the parameter of unit storage-based time cost, $\beta$ is the weight of delay cost in Equations \ref{8} and \ref{9}.

\begin{equation}
C_{Po} = EP \times \sum_{k=1}^{K} E_k
\label{8}
\end{equation}

\begin{equation}
\small
\hspace{-0.3cm}
C_{time} = f_{C_O}(T(waiting)) + \beta \times f_{C_ {cl}^ i}(T(delay)) 
\label{9}
\end{equation}

Waiting time costs are directly influenced by inventory holding costs~\citep{Gupta202}, which typically comprise approximately rate $ \gamma\simeq25\% $of annual inventory expenses~\citep{Goh2014}.  This relationship is expressed linearly in Equations \ref{e10} and \ref{e11}, where $UIHC$ denotes the unit inventory holding cost.

\begin{equation}
\small
f_{C_O}(T(waiting)) = \textit{UIHC} \times T_{wi}  
\label{e10}
\end{equation}

\begin{equation}
\hspace{-0.3cm}
\small
\text{UIHC(min)} = \gamma \times \frac{Average Annual Inventory Cost}{365\times 24\times 60\times Package Quantity } 
\label{e11}
\end{equation}

\begin{figure}[h]
\centering
\subfloat[\footnotesize Expedite.]{\includegraphics[width=1.4in]{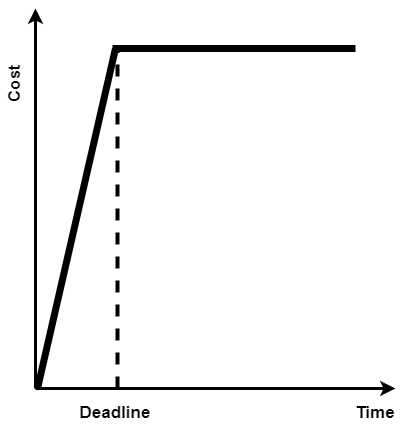}%
\label{fig:sub1}}
\subfloat[\footnotesize Fixed Date.]{\includegraphics[width=1.4in]{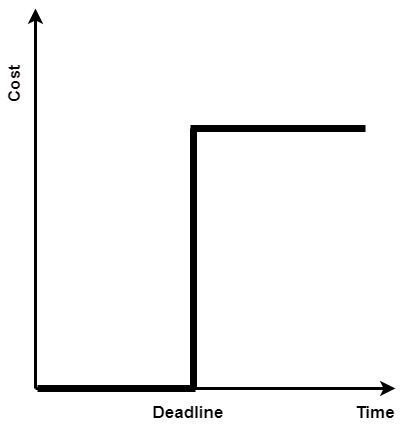}%
\label{fig:sub2}}
\subfloat[\footnotesize Standard Urgency.]{\includegraphics[width=1.4in]{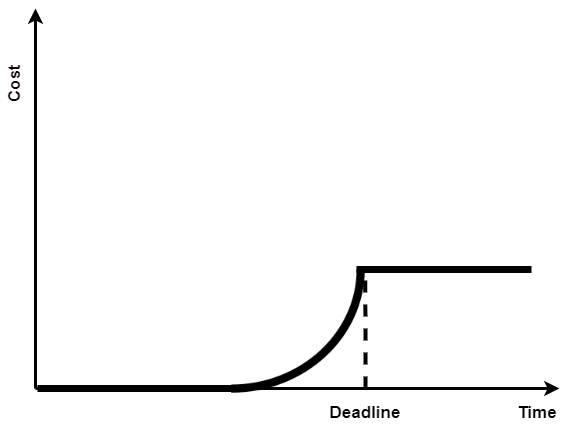}%
\label{fig:sub1}}
\subfloat[\footnotesize Intangible.]{\includegraphics[width=1.4in]{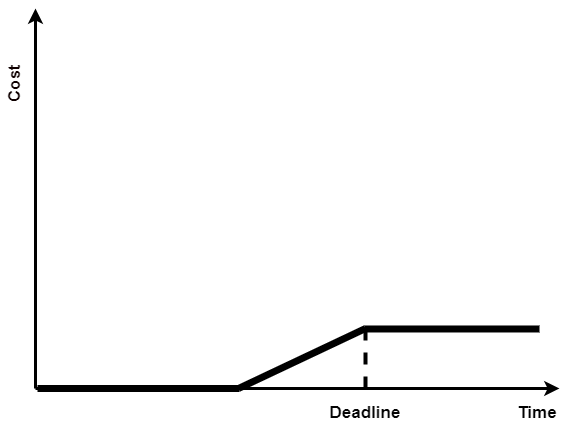}%
\label{fig:sub1}}
\caption{\textbf{ Delay curves of four classes. }}
\label{fig:four delay curves}
\end{figure}

The unit delay cost represents the dollar value of a feature, reflecting the value of the information discovered and how that value decays over time~\citep{Arnold2013}. It is calculated by assessing the impact of the potential losses. In this study, order types are categorized into four delay cost profiles: expedite, standard urgency, intangible, and fixed date ~\citep{Anderson2012}. Expedited orders are high-priority and require immediate service, with zero tolerance for delays. Fixed Date orders are time-sensitive but allow flexibility, incurring costs only when deadlines are missed. Standard Urgency orders gradually accumulate costs, which then taper off to a stable level. Intangible orders, which are low-priority or optional, tolerate extended delays and may even be canceled if necessary. Order classification is determined by delivery attributes, urgency, service deadlines and strategic supply chain positioning. Each class has a predefined deadline, rate of cost escalation, and maximum allowable loss, as shown in Figure ~\ref{fig:four delay curves} and formalized in Equation \ref{eq:12}.
\begin{subequations} \label{eq:12}
\begin{equation}
\small
\hspace{-0.3cm}
f_{C_ {cl}^ a}(T(delay)) =
  \begin{cases}
\lambda_1 \times T(delay) & ,\, T_{wi} < T_{ddi} \\
C_A & ,\, T_{wi} \geq T_{ddi} \\
  \end{cases}
\label{eq:12a}
\end{equation}

\begin{equation}
\small
\hspace{-0.3cm}
f_{C_ {cl}^ b}(T(delay)) =
\begin{cases}
0 & ,\, T_{wi} < T_{ddi} \\
C_B & ,\, T_{wi} \geq T_{ddi} \\
\end{cases}
\label{eq:12b}
\end{equation}

\begin{equation}
\small
\hspace{-0.3cm}
f_{C_ {cl}^ c}(T(delay)) =
\begin{cases}
0 & ,\, T_{wi} < T_{ddi} \\
e^{\lambda_2 \times T(delay)} & ,\, T_{ddi} < T_{wi} < T_u \\
C_C & ,\, T_{wi} \geq T_u \\
\end{cases}
\label{eq:12c}
\end{equation}

\begin{equation}
\small
\hspace{-0.3cm}
f_{C_ {cl}^ d}(T(delay)) =
\begin{cases}
0 & ,\, T_{wi} < T_{ddi} \\
\lambda_3 \times T(delay) & ,\, T_{ddi} < T_{wi} < T_u \\
C_D & ,\, T_{wi} \geq T_u \\
\end{cases}
\label{eq:12d}
\end{equation}
\end{subequations}

\subsubsection{Constraints}\label{constraints} 
The constraint (\ref{flow_conservation}) keeps flow conservation, constraint (\ref{assignment}) demonstrates order assignment, constraint (\ref{subtour}) is a flow conservation constraint which is fundamental to ensuring that AGV routes are continuous and valid circuits, constraint (\ref{capacity}) show the limitation of AGV capacity, constraint (\ref{G}) illustrates the load of orders is under AGV's weight limitation, constraint (\ref{tolerance}) defines total processing duration fitting within the fixed time window, constraint (\ref{time}) makes sure the deadline of orders are later than arriving time, constraint (\ref{pri}) explains the sequence of different classes for deadlines,  constraint (\ref{w}) constrains the weighting factor between energy consumption and service time to remain within the normalized range of 0 to 1, and constraint (\ref{class}) assigns stricter delivery deadlines and higher delay cost penalties to higher priority orders, reflecting their urgency.
\begin{equation}
\sum_{j \in \mathcal{N}} x_{ijk} - \sum_{j \in \mathcal{N}} x_{jik} = 0, \quad \forall i \in \mathcal{N}, \forall k \in \mathcal{K}
\label{flow_conservation}
\end{equation}

\begin{equation}
\sum_{k=1}^{K} y_{ik} = 1, \quad \forall i \in \mathcal{O}
\label{assignment}
\end{equation}

\begin{equation}
u_{ik} - u_{jk} + n \times x_{ijk} \leq n - 1, \quad \forall i,j \in \mathcal{O}, i \neq j, \forall k \in \mathcal{K}
\label{subtour}
\end{equation}

\begin{equation}
\sum_{i=1}^{n} g_i \cdot y_{ik} \leq M, \quad \forall k \in \mathcal{K}
\label{capacity}
\end{equation}

\begin{equation}
\small
\sum_{i=1}^{m} g_i \leq G_L, \forall i \in P_i
 \label{G}
\end{equation}

\begin{equation}
\sum_{k=1}^{K} \sum_{(i,j) \in A} t_{ij} \times x_{ijk} \leq T_B
\label{tolerance}
\end{equation}

\begin{equation}
\small
t_{di} \geq t_{ai}, \forall i \in P_i
 \label{time}
\end{equation}

\begin{equation}
T_{ddi} \leq T_{ddj} \quad \forall i \in C_a, j \in C_b \cup C_c \cup C_d
\label{pri}
\end{equation}

\begin{equation}
\small
0 \leq w_t \leq 1
 \label{w}
\end{equation}

\begin{equation}
\small
C_A \geq C_B \geq C_C \geq C_D
 \label{class}
\end{equation}

\subsection{Priority Rules}\label{Rule} 
\begin{figure*}[htbp]
\centering

\subfloat[\footnotesize FCFS]{%
\includegraphics[width=2.4in,height=1.8in]{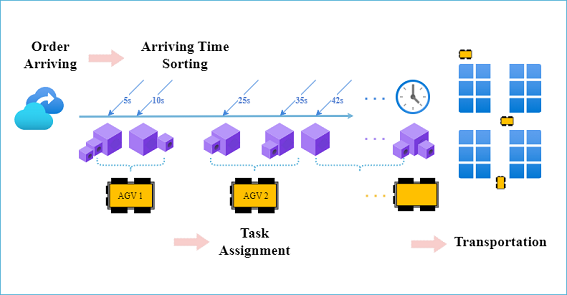}
\label{fig:fcfs}}
\hspace{0.2cm}
\subfloat[\footnotesize SPT]{%
\includegraphics[width=2.4in,height=1.8in]{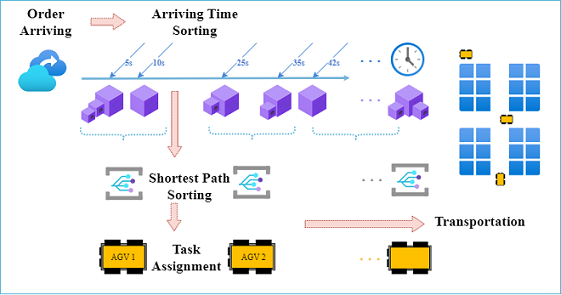}
\label{fig:spt}} \\[0.3cm]

\subfloat[\footnotesize EDT]{%
\includegraphics[width=2.4in,height=1.8in]{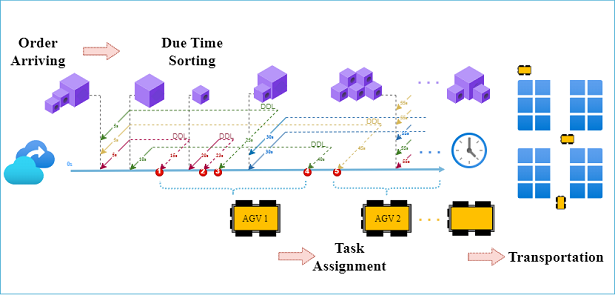}
\label{fig:edt}}
\hspace{0.2cm}
\subfloat[\footnotesize LDC]{%
\includegraphics[width=2.4in,height=1.8in]{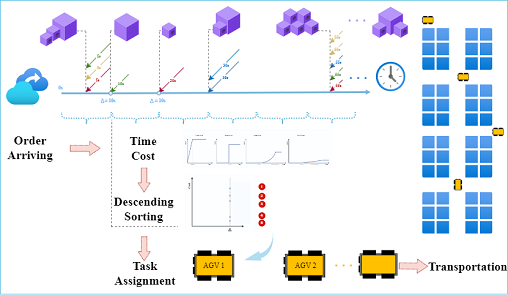}
\label{fig:ldc}} \\[0.3cm]

\subfloat[\footnotesize PDSP]{%
\includegraphics[width=2.4in,height=2in]{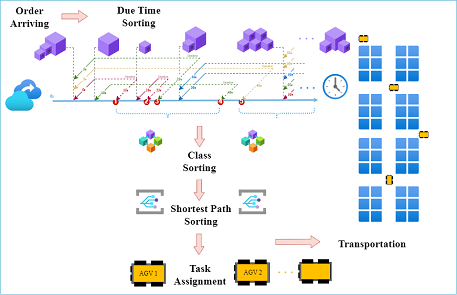}
\label{fig:pdsp}}
\hspace{0.2cm}
\subfloat[\footnotesize DCSP]{%
\includegraphics[width=2.4in,height=2in]{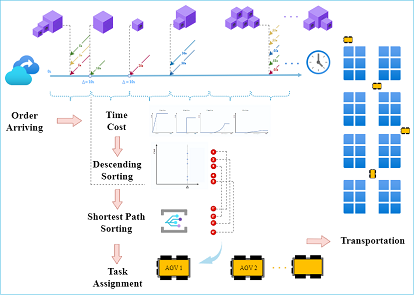}
\label{fig:dcsp}}

\caption{\textbf{ Sorting rules.}}
\label{fig:six_rules}
\end{figure*}

This study incorporates four classical sorting rules from operations research and introduces two new strategies tailored for AGV-based warehouse logistics. The FCFS rule \citep{alma992404284905151} schedules orders strictly by arrival time, assigning each AGV four sequential tasks per cycle. The SPT rule~\citep{Smith1956} schedules based on arrival time and the shortest travel distance, minimizing the path length in each cycle. The Earliest Due Time (EDT) rule~\citep{Goldberg1977} prioritizes tasks with the nearest deadlines, aiming to reduce tardiness, although its correlation with the overall service cost is limited. The Least Delay Cost (LDC) rule assigns tasks based on projected delay penalties, thereby reducing the financial impact of service delays. Based on these, the PDSP rule integrates external priorities, such as customer class and deadlines, with internal route optimization to minimize the travel distance. In contrast, the DCSP rule omits customer classification and prioritizes tasks solely based on the delay cost before optimizing the route sequence. Both proposed methods improve the service quality by balancing the energy efficiency, inventory cost, and delivery performance. A detailed flowchart of the six rules is presented in Figure ~\ref{fig:six_rules}.

\subsection{Simulation-Based Optimization}\label{Rule} 
\subsubsection{Frameworks of Simulation}\label{simulation}
We developed a simulation-based optimization framework using A*guiding multi-agent deep Q-network (A*guiding MADQN) pathfinding and six sorting heuristics within a Python-based transport simulator, extending by \cite{Luo2023} with dynamic priorities and capacity constraints. The simulator processes dynamically generated orders with coordinates, time windows, weight, and priority: {$ (x_i, y_i, t_{oi},t_{ddi},g_i,c_i )$}. These nodes are fed into A* navigation for local delivery routing, and the data are trained using DQN. Nonlinear delay-cost functions model real-time penalties in four priority classes. Each order’s penalty parameters are updated periodically or upon new arrivals for dynamic service evaluation. Two task-sorting strategies are proposed.\\
\textbf{PDSP-AGMADQN:}
This method integrates external priorities (customer classes) and internal routing efficiency. Orders are first sorted by class (A to D) and then by a tolerance window (deadline minus arrival time). After grouping tasks by AGV capacity, A* computes a minimal distance route by optimizing node adjacency and avoiding redundant travel. The data were trained using the MADQN method to achieve the simultaneous collaboration of multiple AGVs.\\
\textbf{DCSP-AGMADQN:}
This approach ignores the customer class and prioritizes tasks purely based on delay cost, dynamically sorting orders based on nonlinear penalty functions. A* then determines the shortest route for the selected tasks and trains the model again. Owing to DCSP's sensitivity of the DCSP to time lags, careful scheduling is required to avoid costly deadline violations.
\subsubsection{Algorithm Application}\label{simulation}
Enhanced A* variants for AGV applications include variants such as dynamic window integration \citep{Chen2022}, hybridized models with ant colony optimization (ACO), the dynamic window approach (DWA) by \cite{Li2023} and congestion-aware priority extensions by \cite{Zhang2020}.Combining this approach with reinforcement learning algorithms has also yielded promising results. In terms of convergence, reward calculation, and loss function performance, the AGMADQN algorithm \citep{Luo2023} used in this study is one of the best to date, to the best of our knowledge. The main idea of the algorithm is to use A* to generate the shortest path and apply the trained neural network to a test structure. By combining two convolutional layers and two linear computations, it extracts features from the warehouse environment (obstacles, shelves, other AGVs, etc.), the position and motion information of the AGV (position coordinates, movement direction, etc.). The adaptability of the method was verified by varying the conditions. This method requires only the current system layout information to guide the AGV through a series of randomly assigned tasks, significantly reducing the training time and increasing the accumulated rewards, thereby ensuring algorithm stability and accuracy. A* estimates the total cost $f(n)$ as:
\begin{equation}
\small
f(n) =  g(n) + h(n)
\label{A*}
\end{equation}
where $g(n)$ is the actual cost from the start node to the node $n$, $h(n)$ means the heuristic estimated cost from n to the destination node.

The state representation function for AGV shows more information and the state $\mathbf{s}_i^t$ for AGV $i$ at time $t$ is defined as:
\begin{equation}
\mathbf{s}_i^t = \Big( \mathbf{p}_i^t, \mathbf{v}_i^t, \mathcal{B}, \mathcal{S}^t, \mathcal{G}_i^t, \mathcal{T}_i^t, \mathcal{A}_{-i}^t \Big)
\label{state function}
\end{equation}
where $\mathbf{p}_i^t=(x_i^t, y_i^t) \quad$ represents the current position and $ \mathbf{v}_i^t=(v_x^t, v_y^t) \quad$ shows the velocity of AGVs. For environment, we set obstacles as static barriers and shelves in $\mathcal{B}=\{(x_k, y_k)\}_{k=1}^K \quad$ and ${S}^t=\{(x_l^t, y_l^t)\}_{l=1}^L \quad $. 
Furthermore, AGVs follow multiple trips for several tasks and finally move to the picking location. The goal hierarchy is shown in Equation \ref{goals}, which illustrates the sub-goals for multiple orders and the final goal position of the picking station at the end of the journal.
\begin{equation}
    \mathcal{G}_i^t = \Big( \underbrace{\{\mathbf{g}_1, \mathbf{g}_2, \mathbf{g}_3, \mathbf{g}_4\}}_{\text{Sub-goals}}, \underbrace{\mathbf{g}_{\text{pick}}}_{\text{Final goal}} \Big)
\label{goals}
\end{equation}
As for task attributes, we define the labels for each AGV in Equation \ref{order} and collect other AGVs' states in $\mathcal{A}_{-i}^t = \{\mathbf{p}_j^t, \mathbf{v}_j^t\}_{j \neq i} \quad$.

\begin{equation}
    \mathcal{T}_i^t = \big( \text{priority}, \text{deadline}, \text{cost}, \text{visited-subgoals} \big)
\label{order}
\end{equation}

The deep q-network (DQN) training iteration follows the update rule in Equation \ref{qiteration}, and the loss is calculated using Equation \ref{loss}, where $y_t$ is the target Q-value, $r_t$ is the immediate reward received after taking action $a_t$ state $s_t$, $\gamma$ is the discount factor, and $\theta^-$ is the target network.

\begin{equation}
\small
\hspace{-0.3cm}
y_t =
  \begin{cases}
r_t & \text{if } \mathbf{s}_{t+1} \text{ is terminal}, \\
r_t + \gamma \max_{a'} Q_{\theta^-}(\mathbf{s}_{t+1}, a') & \text{otherwise}.\\
  \end{cases}
\label{qiteration}
\end{equation}

\begin{equation}
L({\theta}) = r_t + \gamma \cdot \max_{a_{t+1}} Q(s_{t+1}, a_{t+1}) - Q(s_t, a_t)
\label{loss}
\end{equation}
The detailed pseudocode for PDSP-AGMADQN and DCSP-AGMADQN is provided in Algorithm \ref{alg:PDSPDCSP}.

\begin{algorithm}
\footnotesize
\caption{PDSP-AGDQN \& DCSP-AGDQN} \label{alg:PDSPDCSP}
\begin{algorithmic}[1]
\Require
    \State Order requests $Q = \{(x_i,y_i,t_{oi},t_{ddi},g_i,c_i)\}_{i=1}^n$
    \State AGV fleet $\{A_k\}$ with capacity and weight limitation
    \State Warehouse environment with obstacles and stations
    \State Heuristic function $f(n) = g(n) + h(n)$
    \State Number of tasks per AGV: $m$

\State \textbf{Initialization:}
    \State Create DQN $Q(s,a;\theta)$ and target network $\hat{Q}(s,a;\theta^-)$
    \State Initialize replay memory $D$ with capacity $N$
    \State Set $\epsilon \gets $, $\xi \gets $ \Comment{Exploration and A* probabilities}

\For{episode $=1$ to $M$}
    \State Reset warehouse environment, order attributes, get initial state $s_0$
    
    \For{$t=1$ to $T_{max}$}
        \State \textbf{Task Assignment:}
        \If{PDSP}
            \State Sort $Q$ by: $c_i$ (priority), then $t_{ddi}$ (deadline)
        \ElsIf{DCSP}
            \State Compute for each task: $cost_i \gets f(cost)$
            \State Sort $Q$ by $cost_i$ 
        \EndIf
        
        \State \textbf{AGV Task Allocation:}
        \For{each $A_k$ in AGV fleet}
            \State $A_k.tasks \gets$ Next $m$ available tasks from sorted $Q$
            \State Mark these tasks as assigned
        \EndFor
        
        \State \textbf{Path Optimization:}
        \For{each $A_k$ in AGV fleet}
            \If{rand() $< \xi$} \Comment{A* guidance}
                \State $path_k \gets$ Optimal path visiting all $m$ tasks using A*
                \State $a_k \gets ConvertPathToAction(path_k)$
            \Else \Comment{DQN selection}
                \If{rand() $< \epsilon$}
                    \State $a_k \gets$ random valid action sequence
                \Else
                    \State $a_k \gets \text{argmax} _{a} Q(s_t,a;\theta)$
                \EndIf    
            \EndIf
            \State $a_t \gets a_t \cup \{a_k\}$
        \EndFor
        
        \State Execute joint action $a_t$, observe $r_t$, $s_{t+1}$, $done$
        \State Store $(s_t,a_t,r_t,s_{t+1},done)$ in $D$
        
        \State \textbf{Learning:}
        \If{$|D| \geq batch\_size$}
            \State Sample batch from $D$
            \State Compute targets: $y_j = r_j + \gamma \cdot \max_{a'}\hat{Q}(s_{j+1},a';\theta^-) \cdot (1-done_j)$
            \State Update $\theta$ via SGD on $\frac{1}{B}\sum(y_j - Q(s_j,a_j;\theta))^2$
            \State Every $C$ steps: $\theta^- \gets \theta$
        \EndIf
        
        \State $\epsilon \gets \max(\epsilon_{min}, \epsilon \cdot \epsilon_{decay})$
        
        \If{$done$} \textbf{break} \EndIf
    \EndFor
\EndFor

    \Return Optimized policy $\pi(s) = \text{argmax}_a Q(s,a;\theta)$

\end{algorithmic}
\end{algorithm}

\section{Case Study}\label{case} 

This study utilized the publicly available e-commerce shipping data (EDA) dataset from the Kaggle platform as a representative case study \citep{Yadav2021}. The dataset contains over 10,000 shipment records and provides comprehensive information on customer orders, product specifications, and delivery conditions. In this study, order quantities were determined based on unique customer IDs and their associated product serial numbers. The warehouse environment is categorized into three scales(small, medium, and large) based on storage capacity. The corresponding order thresholds are adjusted according to the warehouse scale. Spatially, the warehouse layout was partitioned into areas labeled A through F. An analysis of the dataset indicated that Area F had the highest demand density. To minimize the transportation time for high-demand items, Area F is strategically positioned nearest to the conveyor belt. Areas A through D exhibit relatively balanced demand and are therefore uniformly distributed across the remaining grid spaces. Product storage locations are assigned $(x, y)$ coordinates based on the area classification, with multiple items allowed per grid cell. Customer priority is determined using the customer rating attribute of the dataset. Customers rated as level 1 were classified as having the highest external priority (priority A), whereas those with ratings of 4 and 5 were considered to have the lowest external priority (priority D). The delay penalty for priority D customers is modeled using a time-dependent delay cost curve, and the proportion of each customer class adheres to the empirical distribution in the dataset. The product value is used to estimate the inventory holding cost, assuming that delayed deliveries increase storage expenses. The weight of each product (in grams) directly affects the AGV load dynamics. Items are accumulated sequentially until the AGV reaches its full capacity before being dispatched to the conveyor system. Throughout transportation, the energy consumption varies as a function of both the payload and distance traveled. The AGV type employed in this study was a multilevel bin-handling autonomous vehicle with a payload capacity of 270 kg, self-weight of approximately 450 kg, battery capacity of 42 Ah, and voltage of 24 V. Under full load, the AGV can operate for an average of 5 h at a constant speed of 1 m/s \citep{Qiu2015}.

Table \ref{data} provides the summary statistics of the selected data. The dataset does not contain any personal information of customers owing to privacy protection. The following features describe the data.

\begin{table}[h]
\caption{Data Description}\label{data}%
\begin{tabular}{@{}llll@{}}
\toprule
\text{Title} & \text{Description} & \text{Ranges} & \text{Unit}\\
\midrule
ID & Unique identifier of each order & [1,10977] & ---- \\
Warehouse block & Flat 2D warehouse layout of 5 districts & A,B,C,D,E & Categorical \\
Customer rating & Priority classification of customer & [1,5] & Categorical \\
Price of the product & Each price of an order & 96--310 & Dollars \\
Weight of the product & Each weight of an order & 1001--7846 & Grams \\
\botrule
\end{tabular}
\end{table}
\noindent

\begin{figure}[h]
\centering
\subfloat[\footnotesize Small warehouse]{\includegraphics[width=1.2in,height=1.2in]{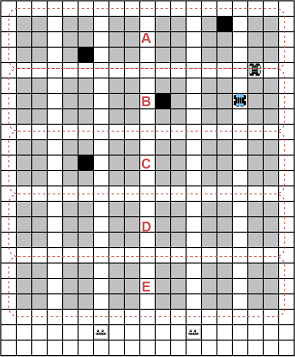}%
\label{fig:sub1}}
\hfil
\subfloat[\footnotesize Middle warehouse]{\includegraphics[width=1.4in,height=1.4in]{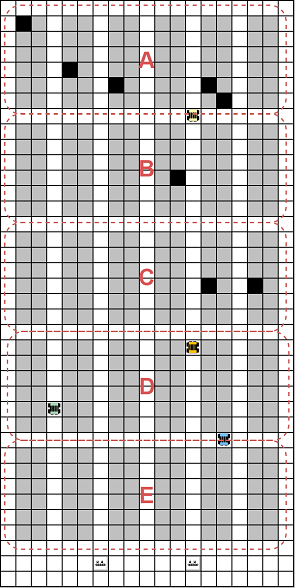}%
\label{fig:sub2}}
\hfil
\subfloat[\footnotesize Large warehouse]{\includegraphics[width=1.6in,height=1.6in]{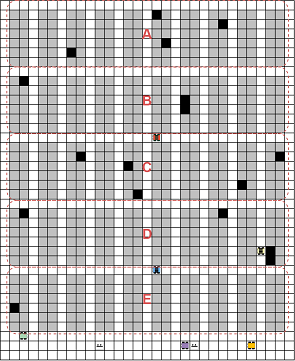}%
\label{fig:sub3}}

\caption{\textbf{Warehouse set-ups.} }
\label{fig:blocks}
\end{figure}

These features collectively inform the simulation model, guiding decisions related to AGV routing, priority scheduling, and cost-sensitive optimization under dynamic warehousing conditions. All the simulation scenarios are shown in Figure ~\ref{fig:blocks}. In these scenarios, the simulation managed congestion when more than one AGV appeared in a close position. The vehicle with a smaller serial number has a higher priority, and other vehicles would stop and wait until the former vehicle completes its last action and leaves.

\section{Results}\label{results}
This section presents a comprehensive evaluation of the proposed methods, including both performance assessment and sensitivity analysis, within the context of the case study.

\subsection{Performance Evaluation}\label{performance}
To evaluate the system performance, 12 key performance indices (KPIs) were defined based on service quality and AGV transportation efficiency. The service metrics included travel time, waiting time, operation time, delay cost, inventory cost, and order-related energy consumption, whereas AGV efficiency was assessed using idle time, running time, and AGV energy use. Two composite indicators, namely the total energy cost and total time cost, captured the overall trade-off, with values recorded in seconds, USD, and Wh. The simulations tested ten scenarios across varying order volumes (1,000 to 10,000) and AGV fleet sizes (1 to 10) using six scheduling methods and class-based delivery penalties that were dynamically updated every 10 s.

\begin{figure}[h]
\centering
\includegraphics[width=0.9\textwidth]{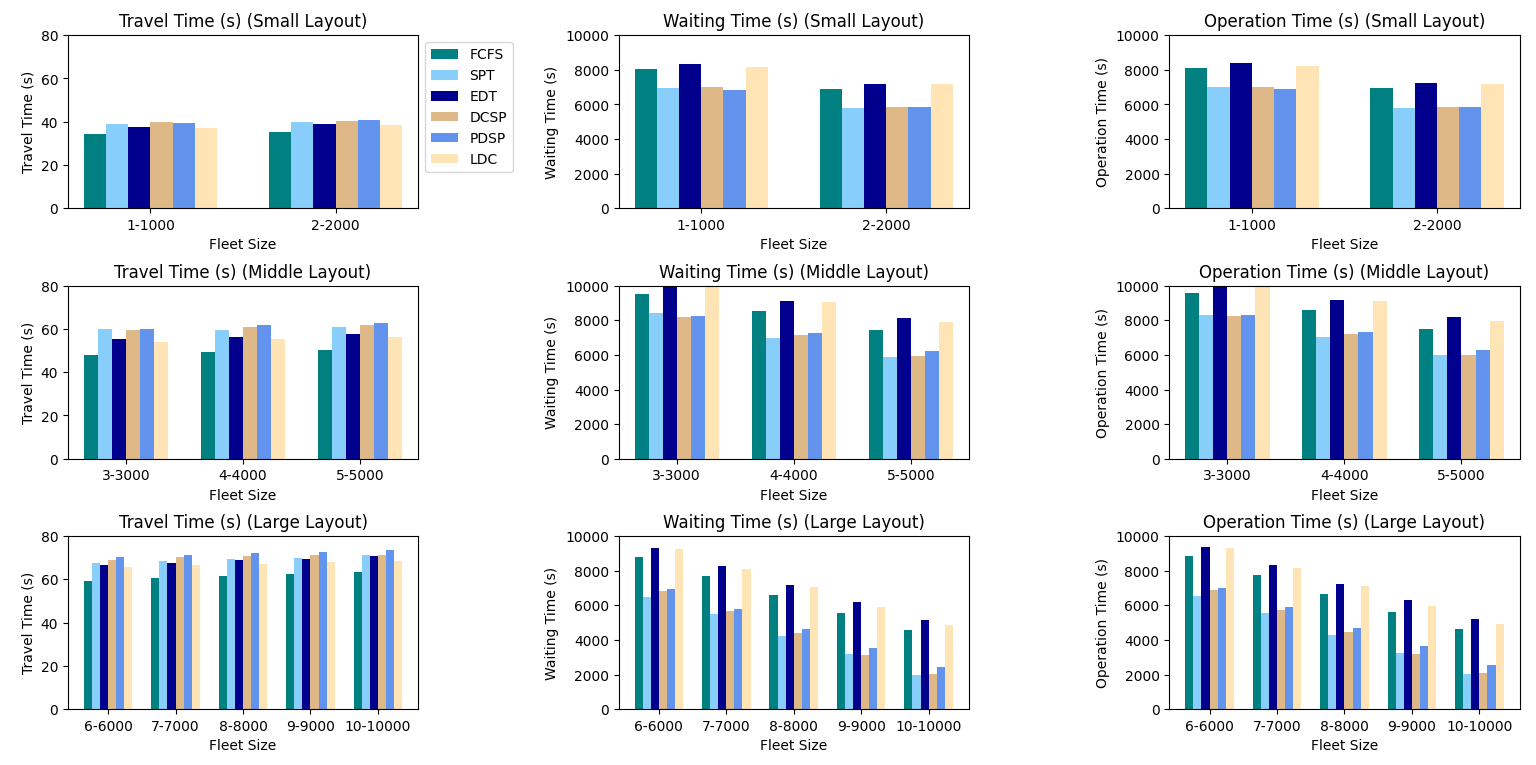}
\caption{\textbf {Service Time.}}
\label{fig:service-time}
\end{figure}

\textbf{Service time} results (cf. Figure \ref{fig:service-time}) show that FCFS achieves the shortest travel time in all scenarios, reducing transport time by 13.34\% on average due to immediate path assignment. However, PDSP and DCSP have the longest travel times because of prioritization and rescheduling. For waiting time, SPT consistently performed the best, whereas EDT exceeded the best methods by 16.46\% to 61.33\%. As fleet size, orders, and warehouse scale increase, the average waiting time decreases significantly by 10\%–20\% in small or medium-scale systems and over 40\% in large-scale systems, indicating higher responsiveness during peak load periods.

\begin{figure}[h]
\centering
\includegraphics[width=0.9\textwidth]{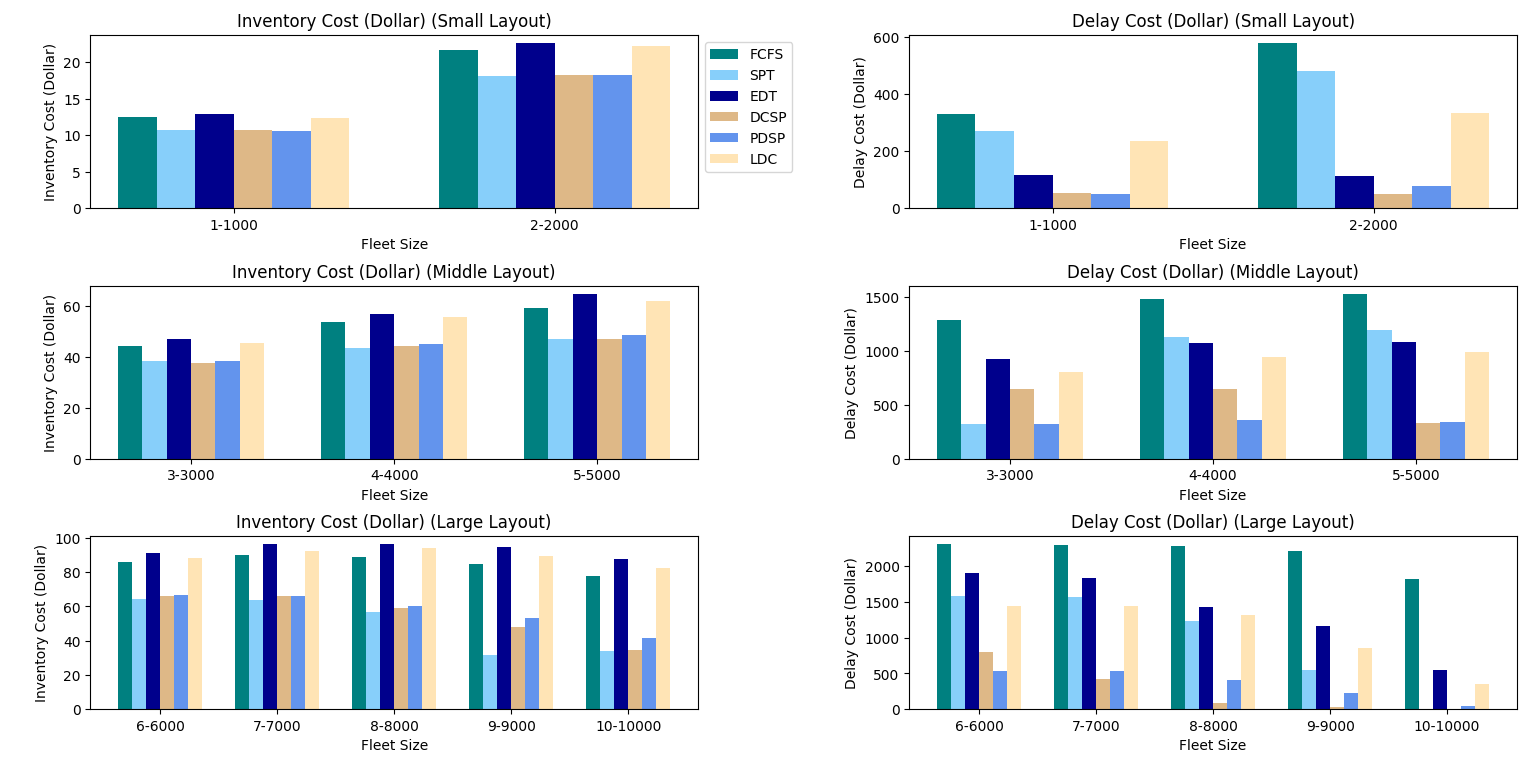}
\caption{\textbf{ Service Cost.}}
\label{fig:service-cost}
\end{figure}

\textbf{Service cost} results (cf. Figure \ref{fig:service-cost}) highlight inventory and delay costs as key economic factors. Inventory costs are mainly driven by waiting time and are highest under EDT and LDC. Delay costs are most sensitive during peak scenarios, with the PDSP and DCSP reducing them by 50\% to 80\%, and up to 90\% in large-scale settings (e.g., 8,000 to 9,000 orders with 8 to 9 AGVs). When the vehicle supply is sufficient (e.g., 10,000 orders with 10 AGVs), both the proposed method and SPT yield minimal delay costs.

\begin{figure}[h]
\centering
\includegraphics[width=0.9\textwidth]{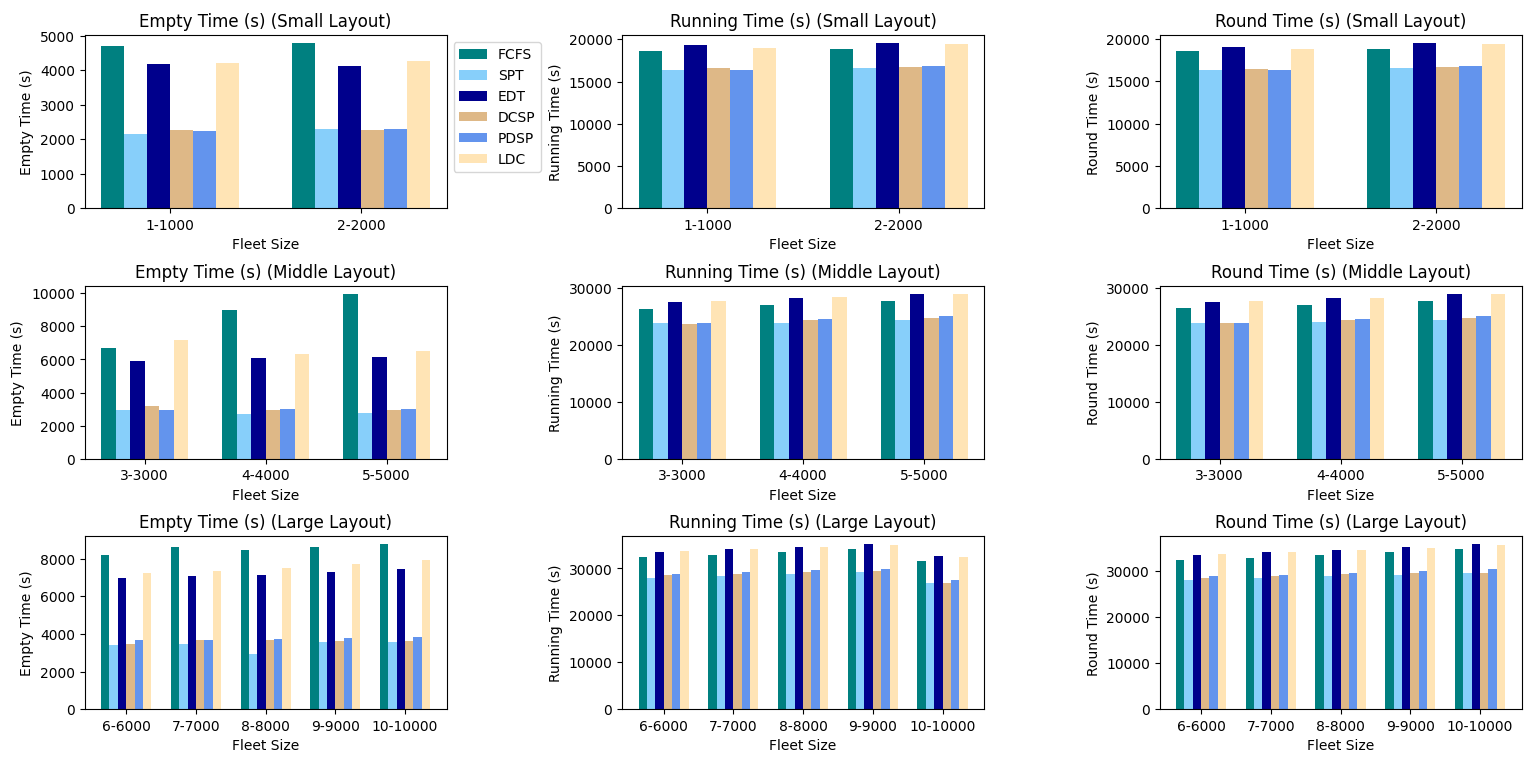}
\caption{\textbf {System Time.}}
\label{fig:AGV&system time}
\end{figure}

\textbf{Transport and system time} findings (cf. Figure \ref{fig:AGV&system time}) reveal FCFS causes the highest AGV idle time due to lack of optimization. In contrast, SPT, DCSP, and PDSP reduce the idle time by 30\% to 50\% via cycle-based order sorting. The AGV running time increases with the number of orders, vehicles, and warehouse scales. SPT achieves the shortest runtime, whereas LDC and EDT achieve the longest.

\begin{figure}[h]
\centering
\includegraphics[width=0.9\textwidth]{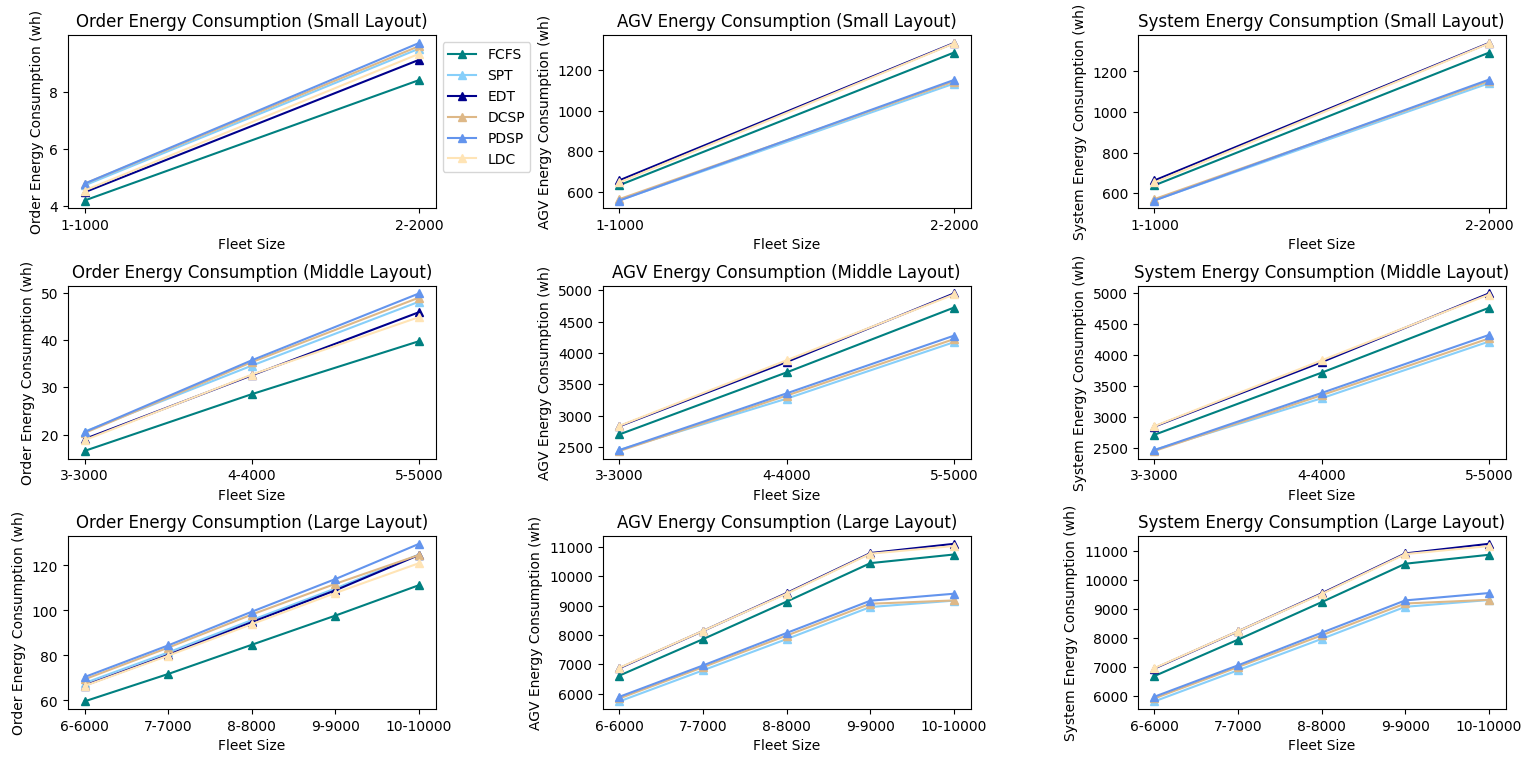}
\caption{\textbf {Energy Consumption.}}
\label{fig:energy}
\end{figure}

\textbf{Energy consumption} analysis (cf. Figure \ref{fig:energy}) separates consumption into goods and AGV movement. FCFS had the lowest order-related energy use, whereas PDSP was often the highest. AGV energy use mirrors the running time trends, with SPT being the most efficient and EDT or LDC the least efficient.

\begin{figure}[h]
\centering
\includegraphics[width=0.9\textwidth]{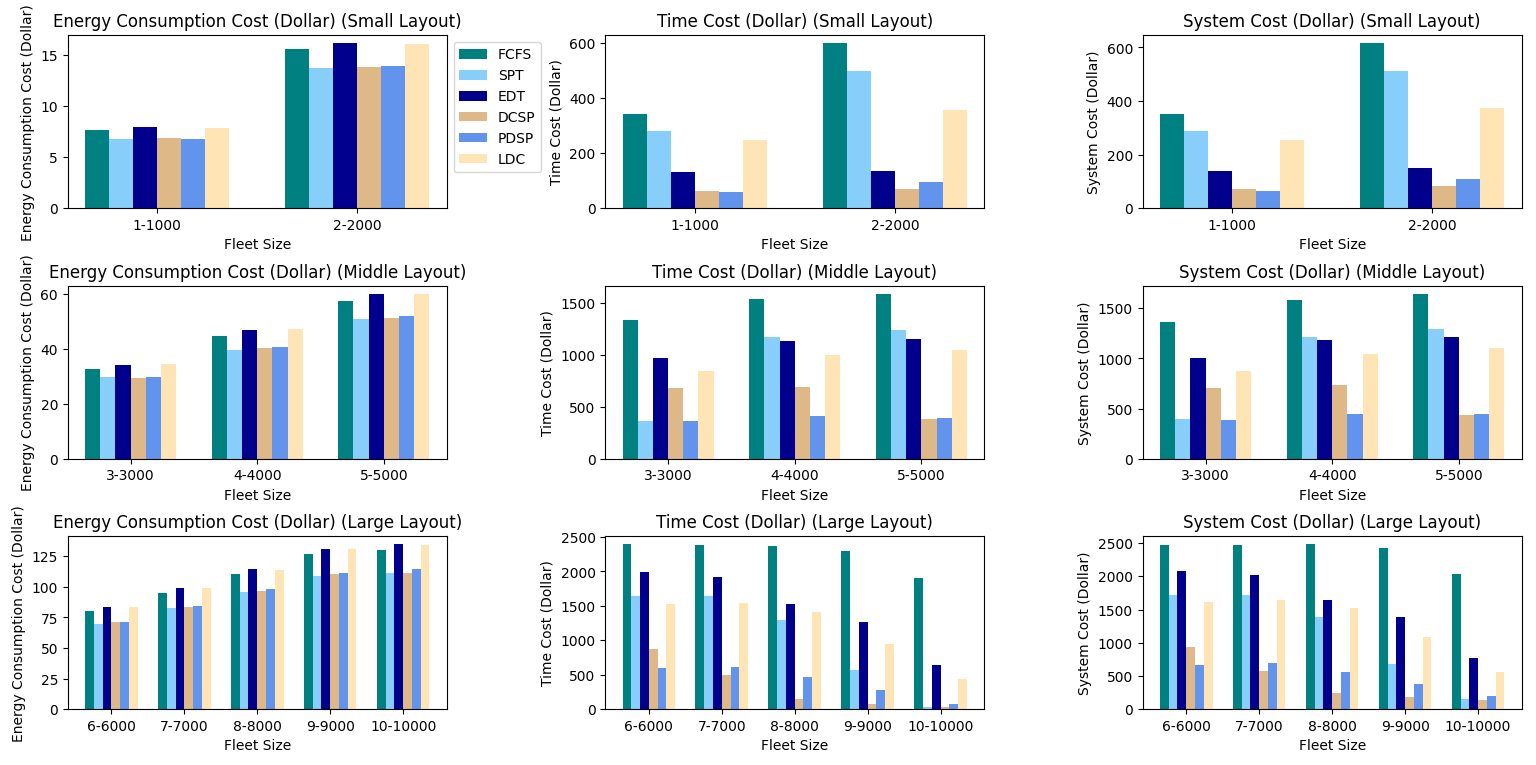}
\caption{\textbf {System Cost.}}
\label{fig:system-cost}
\end{figure}

\textbf{System cost} analysis (cf. Figure \ref{fig:system-cost}) combines time and energy costs equally. The PDSP is the most cost-effective on small or medium scales, whereas the DCSP outperforms as the scale increases, achieving the best results in terms of time cost and second-best in terms of energy cost. Overall, the PDSP or DCSP cut system costs by 50\% to 90\% during high-load periods compared to other methods.

The warehouse scale and layout critically impact these results. As the AGV count and orders increase, larger warehouses increase transport distances and belt loads, driving up time-related costs, such as travel time, waiting time, and delay costs.

\subsection{Sensitivity Analysis}\label{sensitivity} 
Sensitivity comparison between the two methods is shown in Table \ref{DCSP} and Table \ref{PDSP} show how order demand, AGV fleet size, order generation frequency, and delay time windows affect PDSP and DCSP performance across three warehouse scales. Increasing the vehicle count (e.g., from 3 to 5 for 5000 orders or 6 to 10 for 10000 orders) significantly reduces the waiting time, inventory, and delay costs, although the travel time slightly increases (1\% to 2\%) owing to corridor congestion. Increasing the order quantity while holding the fleet size constant (e.g. Five AGVs for 3000 to 5000 orders) leads to proportional increases in all service indicators except travel time, which the PDSP and DCSP methods effectively control. The DCSP cost increase rates (67.56\% and 64.24\%) are higher than those of the PDSP (52.93\% and 44.58\%), showing that the PDSP has a more stable response to rising demand.

\begingroup
\begin{table}[p]
\vspace{0.6cm}
\fontsize{6pt}{6pt}\selectfont
\centering
\begin{sideways}
\begin{minipage}{\textheight}
\caption{Sensitivity Analysis (DCSP)\label{DCSP}}
\setlength{\tabcolsep}{2pt} 
\setlength{\extrarowheight}{4pt} 
\begin{tabular}{p{0.044\textwidth}p{0.014\textwidth}p{0.044\textwidth}p{0.04\textwidth}p{0.08\textwidth}p{0.036\textwidth}p{0.053\textwidth}p{0.053\textwidth}p{0.042\textwidth}p{0.042\textwidth}p{0.048\textwidth}p{0.050\textwidth}p{0.054\textwidth}p{0.054\textwidth}p{0.054\textwidth}p{0.036\textwidth}p{0.048\textwidth}p{0.054\textwidth}p{0.048\textwidth}}
\hline
\multicolumn{5}{c}{Sensitivity parameters}  & \multicolumn{6}{c}{Service} & \multicolumn{3}{c}{Transportation} & \multicolumn{5}{c}{System} \\
Layout & FS & OQ & OWT(s) & DTW(h)  & TrT(s) & WT(s) & OpT(s) & E(1)(wh) &  CoI(\$) & CoD(\$) & EmT(s) & RT(s) & E(2)(wh) & E(wh) & CoE(\$) & CoT(\$) & SRT(s) & CoS(\$)  \\
\hline
\multirow{17}{*}{middle}
        & 3 & 5000 & 0-5 & (1,2,4,4) & 59.74 & 13437.88 & 13497.63 & 47.36 & 103.41 & 1612.65 & 6297.83 & 39733.67 & 4065.23 & 4112.59 & 49.35 & 1716.06 & 39817.50 & 1765.41 \\
        & 4 & 5000 & 0-5 & (1,2,4,4) & 60.73 & 8799.52 & 8860.25 & 48.08 & 68.62 & 1131.56 & 2724.45 & 30381.60 & 4144.54 & 4192.62 & 50.31 & 1200.18 & 30411.67 & 1250.49 \\
        & 5 & 5000 & 0-5 & (1,2,4,4) & 61.86 & 5937.81 & 5999.68 & 48.94 & 47.05 & 334.40 & 2959.22 & 24747.52 & 4219.95 & 4268.89 & 51.23 & 381.45 & 24730.50 & 432.68 \\
\\
        & 5 & 3000 & 0-5 & (1,2,4,4) & 62.57 & 3712.60 & 3775.16 & 21.51 & 17.34 & 12.20 & 1699.80 & 14938.22 & 2547.27 & 2568.78 & 30.83 & 29.54 & 14979.60 & 60.36 \\
        & 5 & 3500 & 0-5 & (1,2,4,4) & 62.50 & 4266.50 & 4329.01 & 27.83 & 23.33 & 34.04 & 2000.38 & 17460.40 & 2977.35 & 3005.18 & 36.06 & 57.37 & 17451.80 & 93.43 \\
        & 5 & 4000 & 0-5 & (1,2,4,4) & 62.10 & 4832.81 & 4894.91 & 35.05 & 30.34 & 61.68 & 2337.04 & 19845.33 & 3384.03 & 3419.08 & 41.03 & 92.02 & 19865.80 & 133.05 \\
        & 5 & 4500 & 0-5 & (1,2,4,4) & 62.05 & 5464.22 & 5526.28 & 42.41 & 38.60 & 99.95 & 2646.24 & 22323.92 & 3806.67 & 3849.08 & 46.19 & 138.56 & 22360.10 & 184.75 \\
\\
        & 5 & 5000 & 0-2 & (1,2,4,4) & 71.00 & 9122.57 & 9184.32 & 48.94 & 71.00 & 1204.56 & 3043.78 & 24765.40 & 4223.00 & 4271.94 & 51.26 & 1275.56 & 24768.60 & 1326.83 \\
        & 5 & 5000 & 0-8 & (1,2,4,4) & 61.52 & 2519.72 & 2581.24 & 48.82 & 20.07 & 6.98 & 2777.98 & 24671.70 & 4207.02 & 4255.83 & 51.07 & 27.05 & 24621.90 & 78.12 \\
\\
        & 5 & 5000 & 0-5 & (1/12,0.3,2,2) & 61.79 & 6057.94 & 6119.73 & 48.89 & 45.23 & 2151.34 & 3134.12 & 24779.92 & 4225.47 & 4274.36 & 51.29 & 2196.58 & 24796.10 & 2247.87 \\
        & 5 & 5000 & 0-5 & (0.3,1,3,3) & 61.92 & 6023.90 & 6085.83 & 49.04 & 46.49 & 1361.73 & 2971.46 & 24685.00 & 4209.29 & 4258.32 & 51.10 & 1408.22 & 24781.00 & 1459.32 \\
        & 5 & 5000 & 0-5 & (1.3,2.3,5.3,5.3) & 61.90 & 5890.08 & 5951.98 & 49.04 & 46.91 & 79.02 & 2912.34 & 24755.10 & 4221.24 & 4270.28 & 51.24 & 125.94 & 24742.30 & 177.18 \\
        & 5 & 5000 & 0-5 & (2,4,8,8) & 62.15 & 5928.39 & 5990.54 & 49.55 & 47.32 & 6.15 & 2787.06 & 24787.64 & 4226.79 & 4276.34 & 51.32 & 53.47 & 24803.00 &104.78 \\
        & 5 & 5000 & 0-5 & (5,10,24,24) & 60.80 & 5847.23 & 5908.03 & 48.09 & 46.66 & 0.00 & 2735.08 & 24426.44 & 4165.20 & 4213.29 & 50.56 & 46.66 & 24445.50 & 97.22 \\
\hline
\multirow{17}{*}{large}
        & 6 & 100000 & 0-5 & (1,2,4,4) & 69.31 & 10947.94 & 11017.25 & 121.66 & 180.28 & 2947.66 & 6746.83 & 47563.23 & 9732.58 & 9854.24 & 118.25 & 3127.95 & 47575.80 & 3246.20 \\
        & 8 & 100000 & 0-5 & (1,2,4,4) & 70.60 & 5345.71 & 5416.31 & 124.42 & 91.08 & 652.78 & 4588.36 & 36473.91 & 9951.25 & 10075.67 & 120.91 & 743.85 & 36476.00 & 864.76 \\
        & 10 & 100000 & 0-5 & (1,2,4,4) & 71.08 & 2039.54 & 2110.61 & 124.68 & 34.61 & 0.08 & 3602.64 & 26912.46 & 9178.23 & 9302.91 & 111.63 & 34.69 & 29628.40 & 146.32 \\
\\
        & 8 & 6000 & 0-5 & (1,2,4,4) & 70.59 & 3419.59 & 3490.19 & 71.29 & 33.49 & 18.84 & 2691.43 & 21854.60 & 5962.63 & 6033.92 & 72.41 & 52.33 & 21896.90 & 124.74 \\
        & 8 & 7000 & 0-5 & (1,2,4,4) & 70.80 & 3859.75 & 3930.55 & 84.47 & 45.04 & 39.75 & 3143.19 & 25546.75 & 6969.97 & 7054.44 & 84.65 & 84.80 & 25588.40 & 169.45 \\
        & 8 & 8000 & 0-5 & (1,2,4,4) & 70.85 & 4379.75 & 4450.59 & 98.20 & 59.00 & 89.24 & 3654.74 & 29251.16 & 7980.65 & 8078.85 & 96.95 & 148.24 & 29292.30 & 245.19 \\
        & 8 & 9000 & 0-5 & (1,2,4,4) & 70.93 & 5002.66 & 5073.59 & 111.48 & 76.42 & 306.71 & 4131.05 & 32944.69 & 8988.37 & 9099.84 & 109.20 & 383.13 & 32990.50 & 492.33 \\
\\
        & 8 & 10000 & 0-2 & (1,2,4,4) & 70.73 & 11804.70 & 11875.43 & 124.30 & 193.52 & 3079.51 & 4731.41 & 36554.76 & 9973.31 & 10097.61 & 121.17 & 3273.03 & 36568.50 & 3394.20 \\
        & 8 & 10000 & 0-8 & (1,2,4,4)& 28.90 & 86.52 & 115.42 & 46.54 & 1.37 & 0.00 & -- & 39582.93 & 10799.49 & 10846.03 & 130.15 & 1.37 & 40436.40 & 131.52 \\
\\
        & 8 & 10000 & 0-5 & (1/12,0.3,2,2) & 71.07 & 5572.47 & 5643.54 & 124.78 & 86.19 & 4384.32 & 4811.00 & 36715.85 & 10017.26 & 10142.04 & 121.70 & 4470.51 & 36700.20 & 4592.22 \\
        & 8 & 10000 & 0-5 & (0.3,1,3,3) & 71.13 & 5559.26 & 5630.39 & 124.76 & 91.69 & 2684.74 & 4739.64 & 36625.76 & 9992.68 & 10117.44 & 121.41 & 2776.43 & 36690.40 & 2897.84 \\
        & 8 & 10000 & 0-5 & (1.3,2.3,5.3,5.3) & 70.49 & 5212.62 & 5283.11 & 99.19 & 71.16 & 99.09 & 4544.34 & 36418.05 & 9936.01 & 10035.20 & 120.42 & 170.24 & 36462.40 & 290.66 \\
        & 8 & 10000 & 0-5 & (2,4,8,8) & 70.10 & 5256.28 & 5326.38 & 98.57 & 71.70 & 11.51 & 4460.40 & 36292.73 & 9901.82 & 10000.39 & 120.00 & 83.21 & 36323.40 & 203.21 \\
        & 8 & 10000 & 0-5 & (5,10,24,24) & 69.12 & 5162.84 & 5231.96 & 97.02 & 70.39 & 0.00 & 4406.85 & 35983.08 & 9817.34 & 9914.35 & 118.97 & 70.39 & 35997.00 & 189.36 \\

\hline
\end{tabular}
\vspace{0.2em}
\begin{minipage}{\textwidth}
\footnotetext{%
\tiny
{Note: FS: fleet size, OQ: order quantity, OTW: time window of orders, DTW: time window of delay, TrT: travel time, WT: waiting time, OpT: operation time; 
CoI: total inventory cost, CoD: cost of order delay, E(1): order energy consumption; 
EmT: AGV empty time, RT: AGV running time, E(2): AGV energy consumption, E: system energy consumption; 
SRT: system round time, CoE: cost of total energy consumption, CoT: cost of total time, CoS: total cost of system.}
}
\end{minipage}
\end{minipage}
\end{sideways}
\end{table}
\endgroup

\begingroup
\begin{table}[p]
\vspace{0.6cm}
\fontsize{6pt}{6pt}\selectfont
\centering
\begin{sideways}
\begin{minipage}{\textheight}
\caption{Sensitivity Analysis (PDSP)\label{PDSP}}
\setlength{\tabcolsep}{2pt} 
\setlength{\extrarowheight}{4pt} 
\begin{tabular}{p{0.044\textwidth}p{0.014\textwidth}p{0.044\textwidth}p{0.04\textwidth}p{0.08\textwidth}p{0.036\textwidth}p{0.053\textwidth}p{0.053\textwidth}p{0.042\textwidth}p{0.042\textwidth}p{0.048\textwidth}p{0.050\textwidth}p{0.054\textwidth}p{0.054\textwidth}p{0.054\textwidth}p{0.036\textwidth}p{0.048\textwidth}p{0.054\textwidth}p{0.048\textwidth}}
\hline
\multicolumn{5}{c}{Sensitivity parameters}  & \multicolumn{6}{c}{Service} & \multicolumn{3}{c}{Transportation} & \multicolumn{5}{c}{System} \\
Layout & FS & OQ & OWT(s) & DTW(h)  & TrT(s) & WT(s) & OpT(s) & E(1)(wh) &  CoI(\$) & CoD(\$) & EmT(s) & RT(s) & E(2)(wh) & E(wh) & CoE(\$) & CoT(\$) & SRT(s) & CoS(\$)  \\
\hline
\multirow{17}{*}{middle}
        & 3 & 5000 & 0-5 & (1,2,4,4) & 60.30 & 13678.51 & 13738.81 & 47.81 & 107.85 & 737.72 & 4900.17 & 40075.47 & 4100.20 & 4148.01 & 49.78 & 845.57 & 40068.70 & 895.35 \\
        & 4 & 5000 & 0-5 & (1,2,4,4) & 61.55 & 8933.21 & 8994.76 & 48.71 & 70.33 & 695.81 & 3728.48 & 30725.45 & 4191.44 & 4240.16 & 50.88 & 766.14 & 30676.40 & 817.02 \\
        & 5 & 5000 & 0-5 & (1,2,4,4) & 62.92 & 6194.13 & 6257.05 & 49.82 & 48.76 & 345.06 & 3033.50 & 25071.54 & 4275.20 & 4325.02 & 51.90 & 393.82 & 25087.60 & 445.72 \\
\\
        & 5 & 3000 & 0-5 & (1,2,4,4) & 62.41 & 3772.02 & 3834.43 & 21.47 & 17.64 & 34.70 & 1840.26 & 14944.14 & 2548.27 & 2569.74 & 30.84 & 52.35 & 14982.40 & 83.19 \\
        & 5 & 3500 & 0-5 & (1,2,4,4) & 62.79 & 4376.55 & 4439.34 & 28.00 & 23.86 & 87.16 & 2110.32 & 17512.16 & 2986.17 & 3014.17 & 36.17 & 111.02 & 17524.50 & 147.19 \\
        & 5 & 4000 & 0-5 & (1,2,4,4) & 62.91 & 4970.25 & 5033.17 & 35.59 & 31.06 & 158.12 & 2404.10 & 20004.56 & 3411.18 & 3446.76 & 41.36 & 189.18 & 20048.70 & 230.54 \\
        & 5 & 4500 & 0-5 & (1,2,4,4) & 62.74 & 5570.27 & 5633.01 & 42.85 & 39.15 & 241.37 & 2754.24 & 22483.30 & 3833.85 & 3876.70 & 46.52 & 280.51 & 22527.60 & 327.03 \\
\\
        & 5 & 5000 & 0-2 & (1,2,4,4) & 62.49 & 9203.11 & 9265.60 & 49.46 & 72.62 & 656.73 & 3047.90 & 24912.30 & 4248.05 & 4297.50 & 51.57 & 729.35 & 24930.90 & 780.92 \\
        & 5 & 5000 & 0-8 & (1,2,4,4) & 63.01 & 2693.65 & 2756.66 & 49.97 & 21.25 & 38.81 & 3033.72 & 25067.20 & 4274.46 & 4324.43 & 51.89 & 60.05 & 25103.90 & 111.94 \\
\\
        & 5 & 5000 & 0-5 & (1/12,0.3,2,2) & 62.62 & 6076.13 & 6138.75 & 49.58 & 47.83 & 550.39 & 3035.54 & 24954.46 & 4255.23 & 4304.82 & 51.66 & 598.22 & 29710.60 & 649.87 \\
        & 5 & 5000 & 0-5 & (0.3,1,3,3) & 62.79 & 6175.85 & 6238.65 & 49.74 & 48.62 & 544.49 & 3047.96 & 25019.96 & 4266.40 & 4316.14 & 51.79 & 593.11 & 25058.30 & 644.90 \\
        & 5 & 5000 & 0-5 & (1.3,2.3,5.3,5.3) & 62.81 & 6130.16 & 6192.96 & 49.75 & 48.27 & 185.06 & 3021.32 & 25034.78 & 4268.93 & 4318.68 & 51.82 & 233.32 & 25064.40 & 285.15 \\
        & 5 & 5000 & 0-5 & (2,4,8,8) & 62.62 & 6148.81 & 6211.42 & 49.64 & 48.40 & 3.50 & 3028.48 & 24989.92 & 4261.28 & 4310.92 & 51.73 & 51.91 & 24993.20 & 103.64 \\
        & 5 & 5000 & 0-5 & (5,10,24,24) & 62.74 & 6138.18 & 6200.92 & 49.67 & 48.33 & 0.00 & 3041.80 & 24987.18 & 4260.81 & 4310.48 & 51.73 & 48.33 & 25019.90 & 100.05 \\
       
\hline
\multirow{17}{*}{large}
        & 6 & 100000 & 0-5 & (1,2,4,4) & 69.93 & 11076.31 & 11146.24 & 122.81 & 183.94 & 836.01 & 6085.68 & 47842.92 & 9789.81 & 9912.61 & 118.95 & 1019.95 & 47844.60 & 1138.90 \\
        & 8 & 100000 & 0-5 & (1,2,4,4) & 71.60 & 5661.96 & 5733.56 & 125.83 & 94.33 & 816.88 & 4690.34 & 36822.28 & 10046.30 & 10172.12 & 122.07 & 911.21 & 36862.80 & 1033.28 \\
        & 10 & 100000 & 0-5 & (1,2,4,4) & 73.72 & 2454.04 & 2527.76 & 129.62 & 41.44 & 40.08 & 3859.81 & 27582.71 & 9406.81 & 9536.42 & 114.44 & 81.52 & 30403.80 & 195.95 \\
\\
        & 8 & 6000 & 0-5 & (1,2,4,4) & 71.97 & 3624.46 & 3696.43 & 72.19 & 34.82 & 130.54 & 2819.14 & 22149.03 & 6042.96 & 6115.15 & 73.38 & 165.37 & 22212.00 & 238.75 \\
        & 8 & 7000 & 0-5 & (1,2,4,4) & 71.98 & 4120.53 & 4192.51 & 85.63 & 46.80 & 250.81 & 3294.10 & 25870.14 & 7058.20 & 7143.83 & 85.73 & 297.62 & 25919.70 & 383.34 \\
        & 8 & 8000 & 0-5 & (1,2,4,4) & 71.94 & 4608.95 & 4680.89 & 99.46 & 60.50 & 404.88 & 3756.85 & 28577.45 & 7796.84 & 7896.31 & 94.76 & 465.37 & 29622.80 & 560.13 \\
        & 8 & 9000 & 0-5 & (1,2,4,4) & 71.53 & 5086.20 & 5157.73 & 112.20 & 75.82 & 582.99 & 4208.29 & 33108.56 & 9033.07 & 9145.28 & 109.74 & 658.81 & 33121.20 & 768.55 \\
\\
        & 8 & 10000 & 0-2 & (1,2,4,4) & 71.62 & 11981.84 & 12053.46 & 125.76 & 199.35 & 1342.67 & 4696.49 & 36848.85 & 10053.55 & 10179.30 & 122.15 & 1542.03 & 36858.90 & 1664.18 \\
        & 8 & 10000 & 0-8 & (1,2,4,4) & 29.68 & 93.35 & 123.03 & 47.16 & 1.48 & 0.00 & -- & 39860.06 & 10875.10 & 10922.26 & 131.07 & 1.48 & 40472.80 & 132.54 \\
\\
        & 8 & 10000 & 0-5 & (1/12,0.3,2,2) & 71.74 & 5663.69 & 5735.43 & 126.02 & 94.38 & 811.35 & 4680.40 & 36871.15 & 10059.63 & 10185.65 & 122.23 & 905.72 & 36913.90 & 1027.95 \\
        & 8 & 10000 & 0-5 & (0.3,1,3,3) & 71.72 & 5616.72 & 5688.45 & 126.11 & 93.62 & 858.00 & 4686.58 & 36904.24 & 10068.66 & 10194.76 & 122.34 & 951.62 & 36907.30 & 1073.96 \\
        & 8 & 10000 & 0-5 & (1.3,2.3,5.3,5.3) & 71.80 & 5668.64 & 5740.43 & 100.96 & 75.57 & 460.75 & 4686.51 & 36925.48 & 10074.45 & 10175.41 & 122.10 & 536.32 & 36936.20 & 658.43 \\
        & 8 & 10000 & 0-5 & (2,4,8,8) & 71.64 & 5598.37 & 5670.01 & 100.63 & 74.66 & 70.11 & 4690.20 & 36831.33 & 10048.77 & 10149.40 & 121.79 & 144.77 & 36873.00 & 266.56 \\
        & 8 & 10000 & 0-5 & (5,10,24,24) & 71.73 & 5602.32 & 5674.04 & 100.78 & 74.70 & 0.00 & 4681.64 & 36810.14 & 10042.98 & 10143.76 & 121.73 & 74.70 & 36853.70 & 196.43 \\
\hline
\end{tabular}
\vspace{0.2em}
\begin{minipage}{\textwidth}
\footnotetext{%
\tiny
{Note: FS: fleet size, OQ: order quantity, OTW: time window of orders, DTW: time window of delay, TrT: travel time, WT: waiting time, OpT: operation time; 
CoI: total inventory cost, CoD: cost of order delay, E(1): order energy consumption; 
EmT: AGV empty time, RT: AGV running time, E(2): AGV energy consumption, E: system energy consumption; 
SRT: system round time, CoE: cost of total energy consumption, CoT: cost of total time, CoS: total cost of system.}
}
\end{minipage}
\end{minipage}
\end{sideways}
\end{table}
\endgroup

When the order frequency increases (e.g., from 0 to 8s to 0 to 2s), the system costs also rise owing to nonstop AGV tasks and higher energy use, whereas slower frequencies allow some idle time but reduce inventory costs, where the DCSP reduces the cost by 77\% and the PDSP by 61\%. Varying delay time windows (from 5 min to 24 h) show that longer windows lower the delay and system costs, with the DCSP showing sharper sensitivity to time pressure. PDSP’s increase is smoother, and both methods offer cost advantages during off-peak or relaxed service windows. Finally, tuning system weights between energy and time (e.g., w = 0.9 for energy focus, w = 0.1 for time-critical scenarios) enhances adaptability, allowing managers to shift priorities based on operational conditions.

\subsection{Algorithms Comparison}\label{sensitivity} 
In this section, we combine the two ranking methods, PDSP and DCSP, with newer advanced reinforcement learning algorithms to compare their adaptability and differences under different algorithms. In this section, we combine the two ranking methods, PDSP and DCSP, with newer advanced reinforcement learning algorithms to compare their adaptability and differences under different algorithms. The main combination algorithms involved are PDSP/DCSP-AGMADQN, PDSP/DCSP-DQN, PDSP/DCSP-AC (Actor Critic), and PDSP/DCSP-A*, which have been used in the literature in this field recently. Considering the computing time and cost, in the algorithm comparison section, we used a small batch of fleets and orders to verify the performance of this method under different reinforcement learning algorithms. The main indicators include operational indicators, such as order completion (throughput per unit time), delay costs, energy consumption costs, and total system costs, as well as algorithm adaptability indicators, such as convergence speed, stability, rewards, and robustness.

\subsubsection{Algorithms Performance}

\begin{figure}[h]
\centering
\subfloat[\footnotesize DCSP-AGMADQN cumulative reward]{\includegraphics[width=3.45in,height=1.6in]{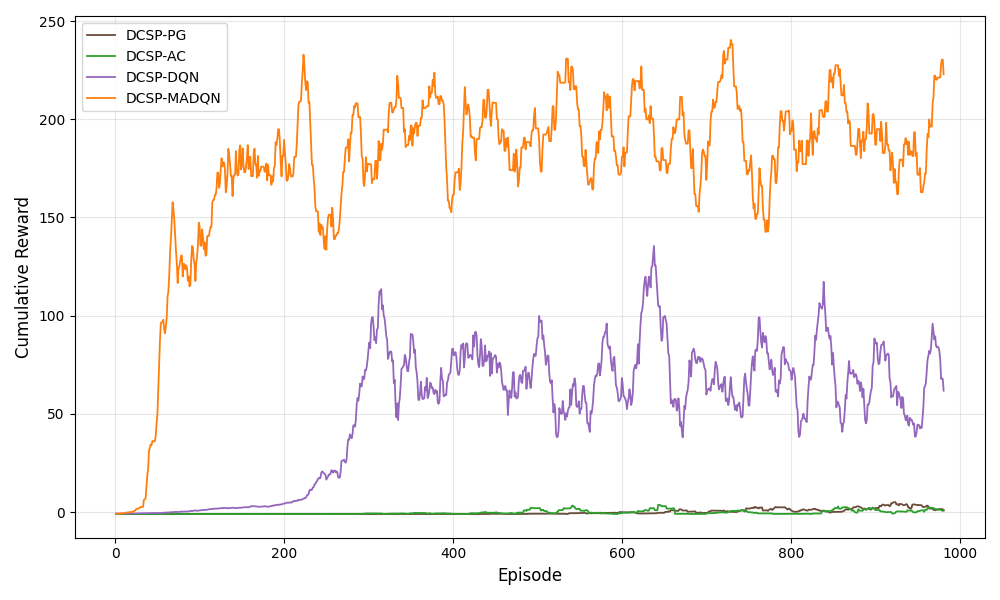}}%
\label{fig:DCSP}
\hfil
\subfloat[\footnotesize  PDSP-AGMADQN cumulative reward]{\includegraphics[width=3.45in,height=1.6in]{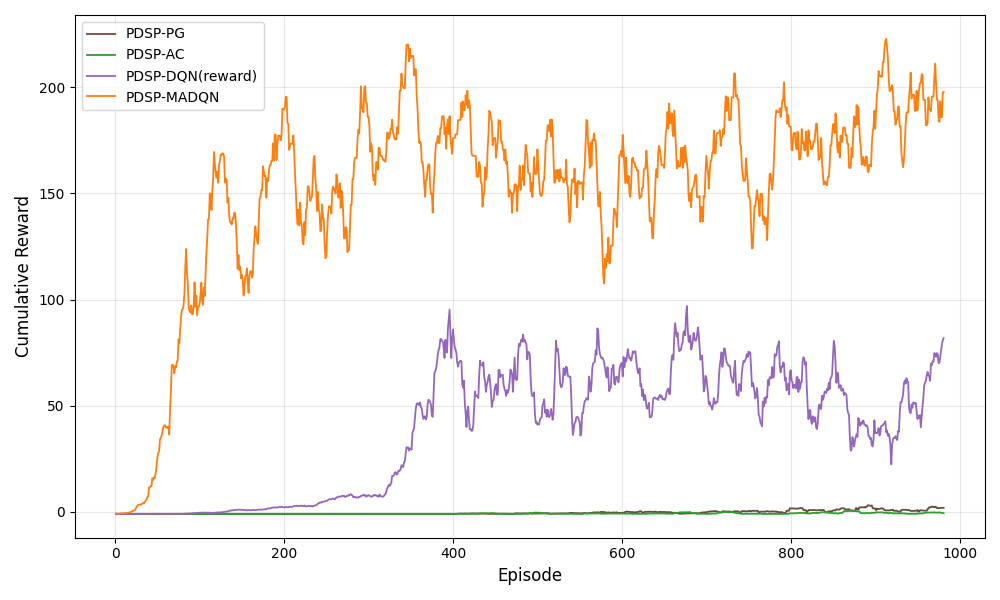}}%
\label{fig:PDSP}
\caption{\textbf{Algorithm Comparison.} }
\label{fig:RL Compare}
\end{figure}

The analysis of cumulative rewards revealed distinct differences in the convergence speed and learning efficiency among the algorithms. In 1000 episodes, we set two AGVs working in a small-scale warehouse for 200 dynamic orders, as shown in Figure \ref{fig:RL Compare}. The standard PDSP/DCSP-DQN exhibited slower and more gradual improvement as it struggled with environmental instability, whereas the PDSP/DCSP-AC algorithm showed moderate speed, balancing the policy and value updates. In contrast, PDSP/DCSP-AGMADQN consistently achieved higher cumulative rewards and maintained a stable growth pattern, indicating its superior capability in handling multi-agent interactions and dynamic scheduling tasks.

The final asymptotic performance, where the cumulative reward curves plateau, indicates the ultimate learning efficiency for 5000 episodes; the algorithm that achieves the highest stable cumulative reward is the most efficient learner. Furthermore, the smoothness of the reward curves is a direct indicator of training stability. A smooth, monotonically increasing curve for AGMADQN suggests high stability.

\begin{figure}[h]
\centering
\includegraphics[width=3.5in,height=1.6in]{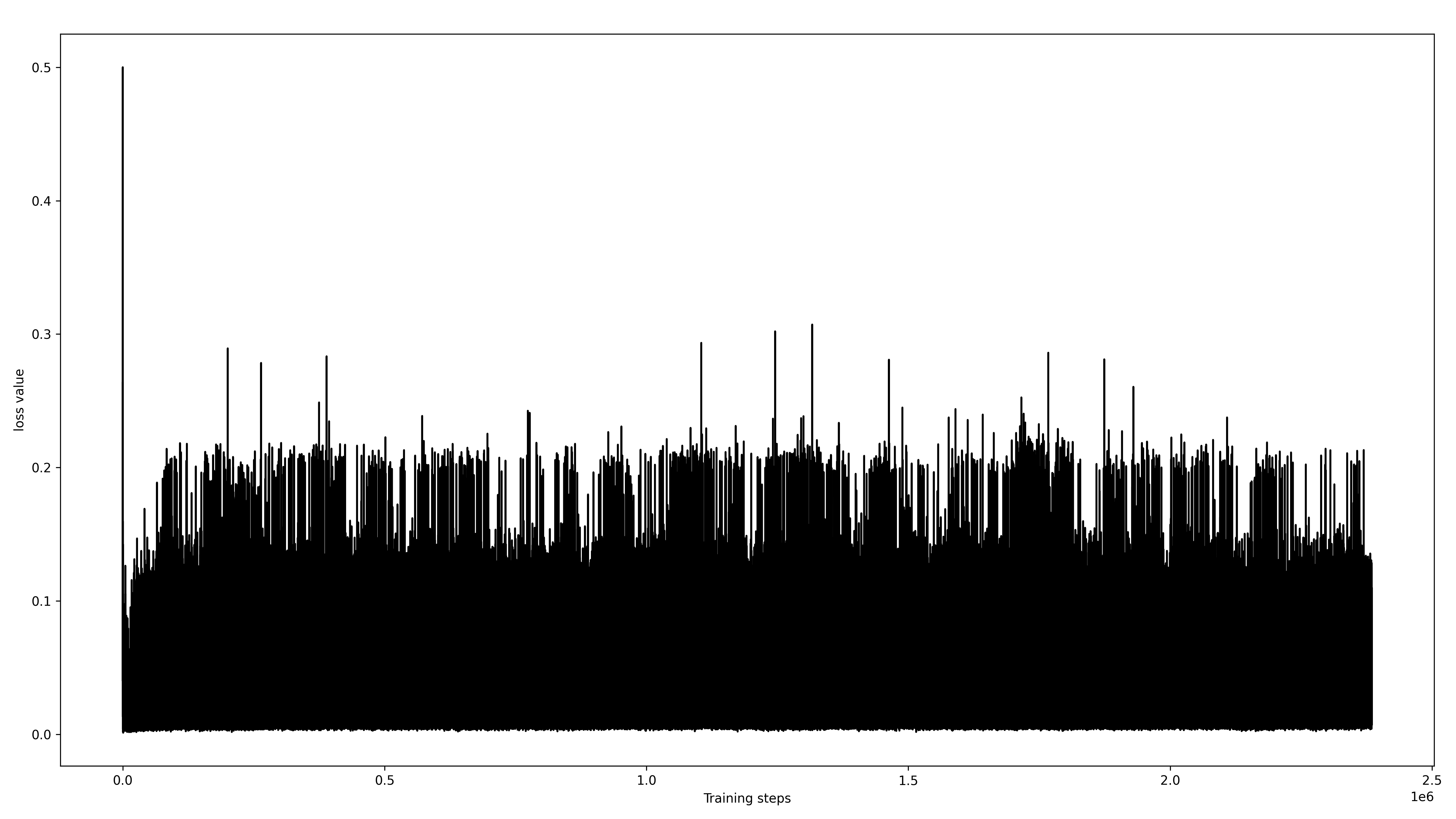}
\caption{\textbf {Loss Value}}
\label{fig:loss}
\end{figure}

Examining the loss values provides insight into the internal learning dynamics and robustness of each algorithm, as shown in Figure \ref{fig:loss}. The rapid decrease and stabilization of the critic (or Q-value) loss in MADQN and AC would correlate with their fast convergence, signifying the quick learning of accurate value estimates. The DQN loss shows greater volatility, reflecting its struggle with a moving data distribution. The robustness of each model can be inferred from its consistent performance across multiple runs with different random seeds. The MADQN, designed for multi-agent settings, should show low variance in both the final reward and loss between runs, indicating robustness to the initial conditions and the inherent randomness of the environment. Conversely, a higher variance in the DQN results highlights its sensitivity and lower robustness. Ultimately, the combination of a high cumulative reward, low final loss, and minimal performance variance positions the MADQN as the superior algorithm in this multi-agent context, effectively balancing speed, stability, and robustness.

\subsubsection{Operational Testing}
\begin{figure}[htbp]
\centering
\subfloat[\footnotesize Testing Order Sequence]{\includegraphics[width=2.6in,height=3in]{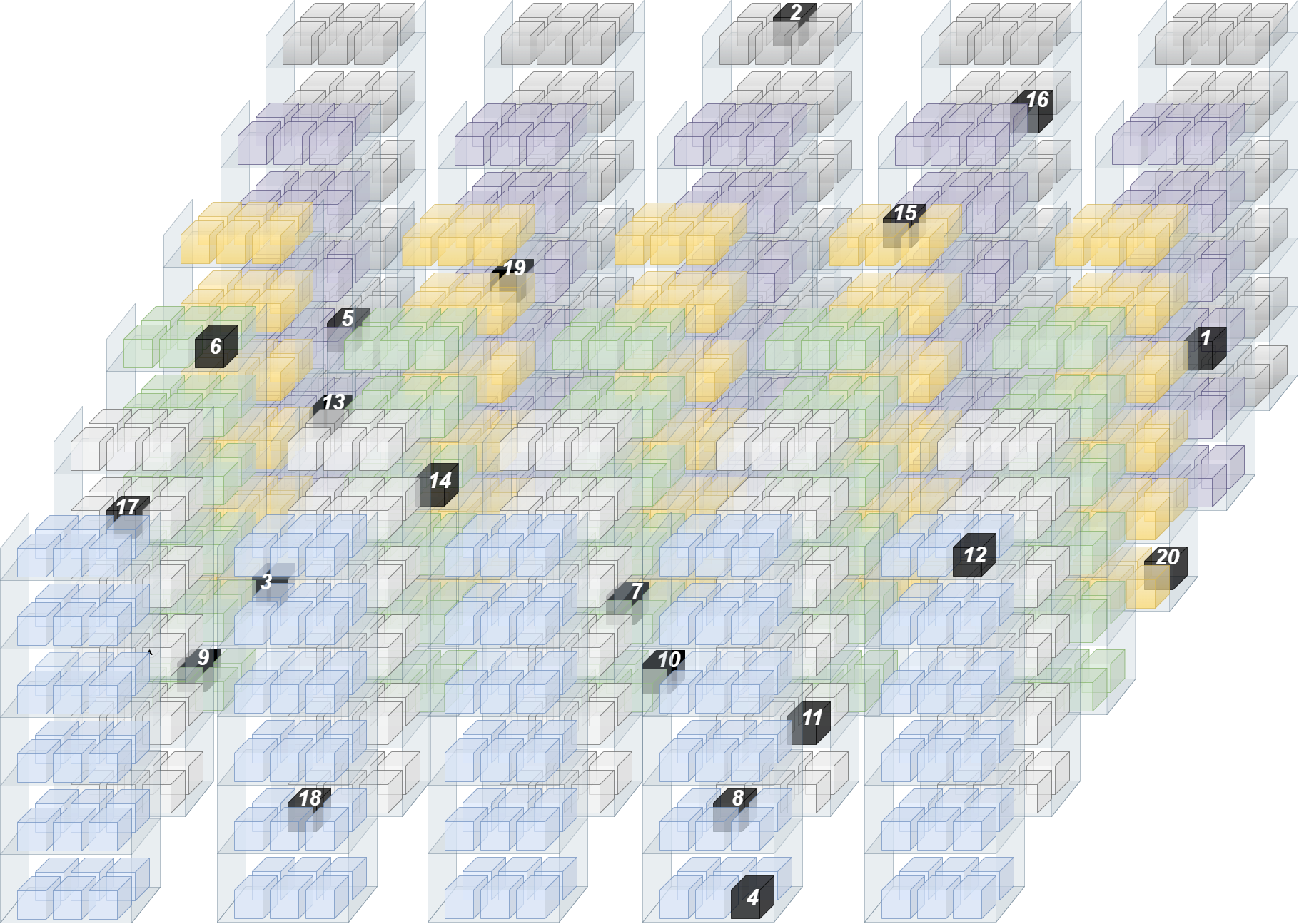}}%
\label{fig:order sequence}
\hfil
\subfloat[\footnotesize  Testing AGV Path]{\includegraphics[width=2.5in,height=3in]{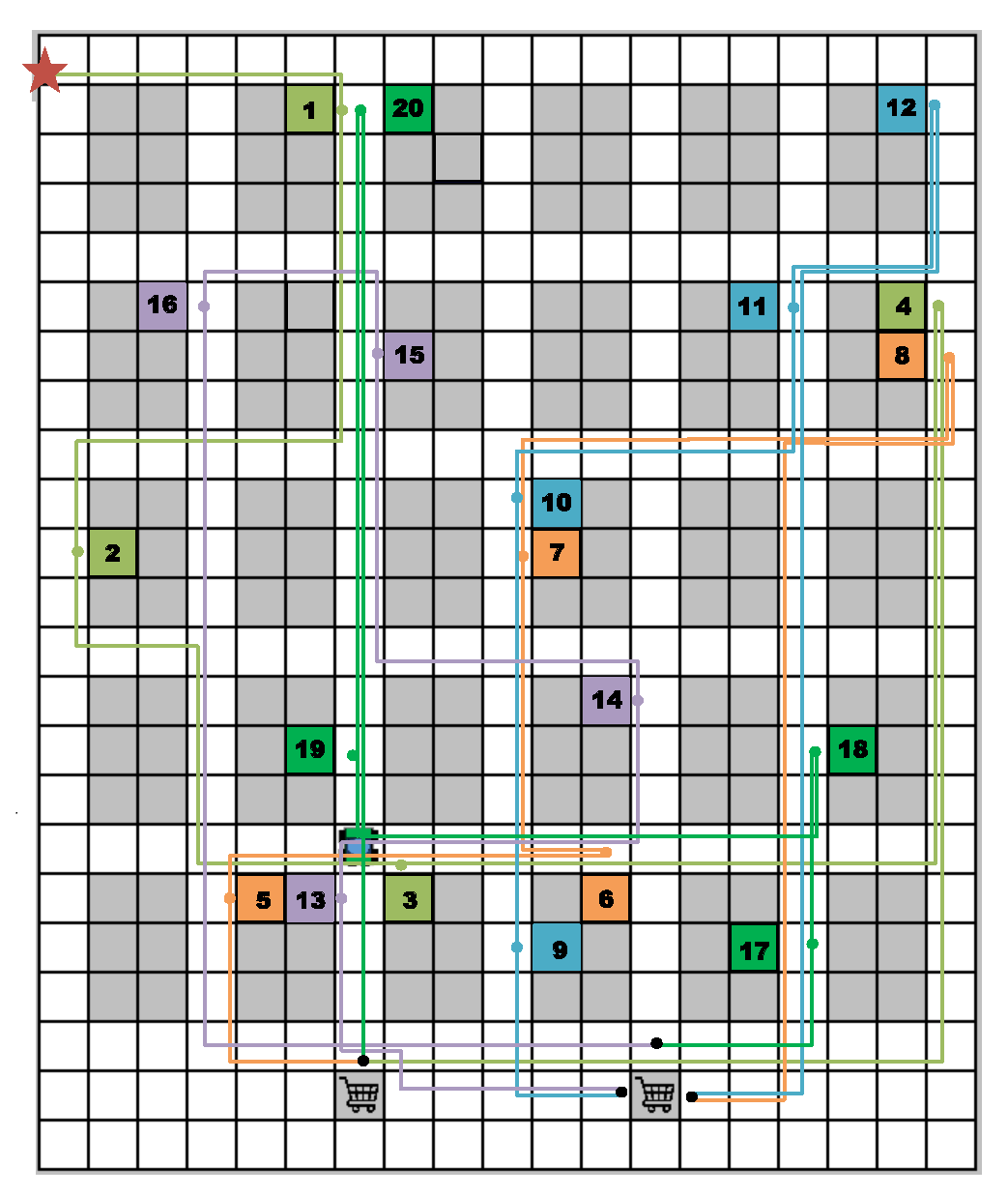}}%
\label{fig:path}
\caption{\textbf{Testing Results.} }
\label{fig:test}
\end{figure}

For operational testing, this paper simulates a scenario where an AGV picks up 20 dynamic orders in a small-scale warehouse. The PDSP-AGDQN algorithm is trained 10,000 episodes and the obtained test results of scheduling and routing policy are compared with those of PDSP-A* as shown in Figure \ref{fig:test}.

The left subplot illustrates the generated order service sequence, demonstrating that the AGV can continuously pick up orders from multiple storage locations according to its training memory and the priority order set by the PDSP. The right subplot visualizes the AGV's corresponding travel path within the warehouse grid environment. After sufficient training, the agent's learned strategy generates a continuous and collision-free path, efficiently connecting all assigned order nodes in the warehouse. This trajectory shows that the AGV can successfully complete all test orders by reasonably executing path consistency and reducing unnecessary detours while maintaining low overall cost. In conclusion, the results confirm that, after sufficient training, the agent can generate practical order fulfillment paths in the test scenario.

\section{Conclusion}\label{conclusion} 

This study introduces a bi-level priority sorting strategy for smart warehouse systems, integrating an enhanced DARP framework to optimize AGV routing, scheduling, and energy efficiency. The mathematical model minimizes the operating time and energy consumption by incorporating external priority, delay time, and shortest-path routing. A real-world case study evaluated the system performance in terms of service satisfaction, resource utilization, and operational efficiency.

Simulation results demonstrate that the proposed PDSP and DCSP strategies consistently outperform traditional benchmarks in high-density, peak-time, and service-sensitive scenarios, achieving cost reductions of 50\% to 90\%.  Sensitivity analysis confirms that the PDSP is more stable and robust for large-scale problems than the DCSP while maintaining less delay and energy costs under varying order volumes and time windows. Their dynamic configurability further enables warehouse managers to adjust fleet size, delay tolerance, and priority weights to achieve diverse operational objectives. 

Therefore, we also applied these two methods to more advanced reinforcement learning algorithms to compare their performance and adaptability in different algorithms. Among them, the PDSP/DCSP-AGMADQN method demonstrated the best performance, with good convergence rate, learning ability, and robustness. In addition, the use of reinforcement learning methods can effectively reduce collisions between AGVs, and the algorithm is faster and more accurate than other methods. Owing to the limitations of the computational time and scale of reinforcement learning, this study only discusses cases in small-scale scenarios; however, it is sufficient to demonstrate the superiority of our method. For currently common large-scale problems, such as large scenarios, a large number of orders, and high-density demand, this study recommends that users model and evaluate within the framework of PDSP/DCSP-A*.

Building on this research, we will expand it in future studies. For example, we will transform the fixed delay curve into a dynamic one influenced by external warehouse factors, enabling the system to truly achieve internal and external linkage and ensuring the effectiveness and timeliness of services; we will combine the actual transportation results of the last mile with the optimization of warehouse scheduling to form a closed loop for the automated system; and we will explore more AI or reinforcement learning methods to make the automated logistics system more intelligent and stable.

\section*{Declarations}

\begin{itemize}
\item Funding:
This research is funded by a grant from the Canada Research Chair program in Disruptive Transportation Technologies and Services (CRC-2021-00480) and NSERC Discovery (RGPIN-2020-04492) fund. 
\item Conflict of interest:
The authors declare no competing interests.
\item Ethics approval and consent to participate:
Not applicable. This study uses only publicly available, anonymized data.
\item Data availability: The primary dataset used in this study is the publicly available from the Kaggle platform. The dataset can be accessed at its permanent repository: \url{https://www.kaggle.com/datasets/prachi13/customer-analytics}. During the analysis, exploratory work was also conducted using the public Kaggle notebook \url{https://www.kaggle.com/code/niteshyadav3103/eda-e-commerce-shipping-data}. All data used is fully anonymized.
\item Consent for publication and Data availability:
This study utilizes a publicly available, fully anonymized dataset sourced from the Kaggle platform \cite{}. The data has been thoroughly de-identified to ensure no individual can be identified. Therefore, no additional ethical approval or participant consent is required for the publication of this research.
\item Materials availability:
This is a computational study that did not generate any unique physical materials. All analyses were performed using publicly available digital data and code, the availability of which is detailed in the “Data availability” and “Code availability” statements above.
\item Code availability:
All custom code and algorithms supporting the findings of this study will be made publicly available upon acceptance of the manuscript. The code repository will be hosted on GitHub. During the peer-review process, the code is available from the corresponding author upon reasonable request.
\item Author contribution:
\textit{Xiaozhu Sun}: Conceptualization, Methodology, Validation, Data Curation, Writing- Original draft preparation, Software. \textit{Bilal Farooq}: Resources, Conceptualization, Methodology, Writing- Reviewing and Editing, Funding acquisition, Supervision.
\end{itemize}

\bibliography{sn-bibliography}


\end{document}